\documentclass[reqno,11pt]{amsart} 
\usepackage{amsmath}
\usepackage{amssymb}
\usepackage{amsthm}

\newcommand\NoBlackBoxes{\global\overfullrule0pt}
\NoBlackBoxes
\theoremstyle{plain}

\begin{document}

\title{Spherical covariance representations}

\author{Sergey G. Bobkov$^{1}$}
\thanks{1) 
School of Mathematics, University of Minnesota, Minneapolis, MN, USA
}

\author{Devraj Duggal$^{2}$}
\thanks{2) School of Mathematics, University of Minnesota, Minneapolis, MN, USA
}


\begin{abstract}
Covariance representations are developed for the uniform distributions
on the Euclidean spheres in terms of spherical gradients and Hessians.
They are applied to derive a number of Sobolev type inequalities and to
recover and refine the concentration of measure phenomenon, including
second order concentration inequalities. A detail account is also given
in the case of the circle, with a short overview of H\"offding's kernels
and covariance identities in the class of periodic functions.
\end{abstract} 
 
\maketitle
\markboth{Sergey G. Bobkov and Devraj Duggal}{Periodic covariance representations}

\def\theequation{\thesection.\arabic{equation}}
\def\E{{\mathbb E}}
\def\R{{\mathbb R}}
\def\C{{\mathbb C}}
\def\P{{\mathbb P}}
\def\Z{{\mathbb Z}}
\def\S{{\mathbb S}}
\def\I{{\mathbb I}}
\def\T{{\mathbb T}}

\def\s{{\mathbb s}}

\def\G{\Gamma}

\def\Ent{{\rm Ent}}
\def\var{{\rm Var}}
\def\Var{{\rm Var}}
\def\cov{{\rm cov}}
\def\V{{\rm V}}

\def\H{{\rm H}}
\def\Im{{\rm Im}}
\def\Tr{{\rm Tr}}
\def\s{{\mathfrak s}}
\def\A{{\mathfrak A}}
\def\m{{\mathfrak m}}

\def\k{{\kappa}}
\def\M{{\cal M}}
\def\Var{{\rm Var}}
\def\Ent{{\rm Ent}}
\def\O{{\rm Osc}_\mu}

\def\ep{\varepsilon}
\def\phi{\varphi}
\def\vp{\varphi}
\def\F{{\cal F}}

\def\be{\begin{equation}}
\def\en{\end{equation}}
\def\bee{\begin{eqnarray*}}
\def\ene{\end{eqnarray*}}

\thispagestyle{empty}

\section{{\bf Introduction}}
\setcounter{equation}{0}

\vskip2mm
\noindent
Let $\gamma_n$ denote the standard Gaussian measure on $\R^n$
with density
$$
\varphi_n(x) = (2\pi)^{-n/2}\,e^{-|x|^2/2}, \quad x \in \R^n.
$$
The space $\R^n$ is equipped with the Euclidean norm
$|\cdot|$ and the inner product $\left<\cdot,\cdot\right>$.

It is well known that the covariance functional with respect to
this measure
$$
\cov_{\gamma_n}(f,g) = \int_{\R^n} f g\,d\gamma_n - 
\int_{\R^n} f\,d\gamma_n \int_{\R^n} g\,d\gamma_n
$$
admits the representation
\be
\cov_{\gamma_n}(f,g) = \int_{\R^n} \int_{\R^n} 
\left<\nabla f(x),\nabla g(y)\right> d\pi_n(x,y)
\en
for a certain probability measure $\pi_n$ on $\R^n \times \R^n$.
Here $f$ and $g$ may be arbitrary smooth functions on $\R^n$ with
$\gamma_n$-square integrable gradients $\nabla f$ and $\nabla g$.

In an equivalent form, this identity was derived in the works
by Houdr\'e and P\'erez-Abreu using an interpolation argument
in an infinite-dimensional setting for functionals of the Wiener process
and by Ledoux by means of the Ornstein-Uhlenbeck semi-group,
cf. \cite{H-P}, \cite{Le1}. The form (1.1) with an explicit description of 
the mixing measure $\pi_n$ was later proposed
in \cite{B-G-H}. It was shown that (1.1) can serve as a convenient tool to recover
a number of Sobolev-type inequalities in Gauss space such as the Poincar\'e-type
inequality and exponential bounds, as well as to refine the classical 
dimension-free concentration in the form of a deviation inequality.

The measure $\pi_n$ in (1.1) is unique. Indeed, extending this relation by 
linearity to complex-valued functions and applying it to
$f(x) = e^{i\left<t,x\right>}$, $g(x) = e^{i\left<s,x\right>}$
with parameters $t,s \in \R^n$, we get an explicit expression for 
the Fourier-Stieltjes transform $\hat \pi_n$ of $\pi_n$, namely
$$
\hat \pi_n(t,s) = e^{-\frac{1}{2}\,(|t|^2 + |s|^2)} 
\frac{1 - e^{-\left<t,s\right>}}{\left<t,s\right>}.
$$
It determines $\pi_n$ in a unique way. Although this measure
is not Gaussian, it has $\gamma_n$ as marginals, which is a crucial property 
for many applications.

What is also interesting, only Gaussian measures may satisfy a representation
such as (1.1) in dimension $n \geq 2$:
If for a probability measure $\mu$ on $\R^n$ there exists an identity
\be
\cov_{\mu}(f,g) = \int_{\R^n} \int_{\R^n} 
\left<\nabla f(x),\nabla g(y)\right> d\pi(x,y)
\en
with some (signed) measure $\pi$ on $\R^n \times \R^n$, then necessarily
$\mu$ must be a Gaussian measure with covariance matrix proportional
to the identity matrix (\cite{B-G-H}, Theorem 2.4). Thus, (1.2) characterizes
the class of Gaussian measures.

Nevertheless, it is natural to wonder whether or not one can obtain (1.2) 
with a proper modification of the Euclidean gradient. This turns out to be 
the case where $\mu$ is a uniform distribution on the hypercube $M = \{-1,1\}^n$, 
that is, the product $n$-th power of the symmetric Bernoulli measure. Then 
we have (1.2) for the discrete gradient of functions on $M$ (\cite{B-G-H}).

In this paper, we consider a similar question for the uniform distribution
$\mu = \sigma_{n-1}$ on the unit sphere
$$
\S^{n-1} = \{x \in \R^n: |x| = 1\}, \quad n \geq 2.
$$
Any smooth function $f$ on $\S^{n-1}$ has a (continuous) spherical gradient 
$\nabla_\S f(x)$ which may be defined for every point $x$ on the sphere as the 
shortest vector $w \in \R^n$ satisfying the Taylor expansion up to the linear term
\be
f(x') = f(x) + \left<w,x'-x\right> + o(|x'-x|), \quad x' \rightarrow x  \
(x,x' \in \S^{n-1}).
\en
As we will see, (1.2) does exist for the spherical gradient.

\vskip5mm
{\bf Theorem 1.1.} {\sl On $\S^{n-1} \times \S^{n-1}$ there 
exists a positive measure $\nu_n$ such that, for all smooth functions 
$f,g$ on $\S^{n-1}$,
\be
{\rm cov}_{\sigma_{n-1}}(f,g) = \int_{\S^{n-1}}\int_{\S^{n-1}} 
\left<\nabla_\S f(x),\nabla_\S g(y)\right> d\nu_n(x,y).
\en
Moreover, $\nu_n = c_n \mu_n$ for a probability measure $\mu_n$ on 
$\S^{n-1} \times \S^{n-1}$ with marginals $\sigma_{n-1}$ and with a constant satisfying
$$
\frac{1}{n-1} < c_n < \frac{1}{n-2} \quad (n \geq 3).
$$
}

\vskip2mm
Here $\nu_n$ has density with respect to the product measure 
$\sigma_{n-1} \otimes \sigma_{n-1}$ of the form 
$\psi_n(x,y) = \psi_n(\left<x,y\right>)$ 
for a certain positive function $\psi_n$ on $[-1,1]$. This ensures that 
the marginals of $\nu_n$ represent (equal) multiples of $\sigma_{n-1}$.

One should however emphasize that $\nu_n$ in (1.4) is not unique,
in contrast with the Gaussian case. Let us explain this in the case of the circle.
Let $\kappa$ be any signed measure on the torus $\S^1 \times \S^1$ supported on
$$
\Delta = \big\{(x,y) \in \S^1 \times \S^1: \left<x,y\right> = 0\big\}.
$$
On this set necessarily $\left<\nabla_\S f(x),\nabla_\S g(y)\right> = 0$, 
since the vector $\nabla_\S f(x)$ is orthogonal to $x$, while $\nabla_\S g(y)$ 
is orthogonal to $y$. Hence, if (1.4) holds true for $\nu$, this relation continues
to hold for the measure $\nu + \kappa$. In addition, $\kappa$ may have 
$\sigma_1$ as marginals, for example, when $\kappa$ is
the image of $\sigma_1$ under the mapping $x \rightarrow (x,e^i x)$ from 
$\S^1$ to $\S^1 \times \S^1$.

Despite the non-uniqueness issue, the measure $\nu_n$ in Theorem 1.1 can be
constructed in rather natural ways leading to the same result. 
The first approach is based on the Gaussian covariance identity (1.1) which 
we apply to 0-homogeneous functions $f$ and $g$. For this aim, we give
an integrable description of the density of $\pi_n$ (Section 2) and then integrate
in polar coordinates in order to reduce the right-hand side of (1.1)
to the form (1.4) and thus to obtain an explicit integrable description for the density 
of $\nu_n$ (Sections 3-4). This will also allow us to analyze an asymptotic behavior of
the function $\psi_n(x,y)$ when the points $x,y$ are close to each other and 
to show that
\be
\psi_n(x,y) \sim \frac{1}{|x - y|^{n-3}} \ \ (n \geq 4), \qquad
\psi_3(x,y) \sim \log\frac{1}{|x - y|},
\en
where the equivalence is understood as two sided bounds with $n$-dependent
factors (Section~5). This relation shows that for the growing dimension $n$ 
the mixing measure $\pi_n$ is almost concentrated on the diagonal $x=y$.

In the case of the circle, which is treated separately 
in Section 6, the density $\psi_2$ is bounded. We will return to this case
in the end of the paper in order to clarify the meaning of (1.1),
connect it with the class of periodic covariance identities on the line, 
and eventually evaluate the density on the torus (addressing the uniqueness
problem as well).

In Sections 7-8 we collect basic definitions and identities related to the
differentiation on the sphere, the spherical Hessian and Laplacian 
$\Delta_\S$. Then we discuss a heat semi-group approach to the covariance functional
with respect to $\sigma_{n-1}$, by involving the semi-group $P_t = e^{t \Delta_\S}$ 
(Section 9). This may be done in analogy 
with M. Ledoux' approach to the Gaussian covariance representation, which leads 
to the following counterpart of Theorem 1.1.

\vskip5mm
{\bf Theorem 1.2.} {\sl For all smooth functions $f,g$ on $\S^{n-1}$
$(n \geq 3)$,
\be
{\rm cov}_{\sigma_{n-1}}(f,g) = \int_0^\infty \int_{\S^{n-1}} 
\left<\nabla_\S f,\nabla_\S P_t g\right> d\sigma_{n-1}\,dt.
\en
}

The relationship between the two representations (1.4) and (1.6) are briefly
discussed in Section 10. Let us note here that (1.6) is no longer true on the circle.

Then we turn to applications of (1.4) including the spherical concentration.
Let us recall that, for any function $f$ on $\S^{n-1}$ with 
Lipschitz semi-norm $\|f\|_{{\rm Lip}} \leq 1$ and with $\sigma_{n-1}$-mean $m$, 
there is a deviation inequality
\be
\sigma_{n-1}\big\{|f - m| \geq r\big\} \leq 2e^{-(n-1)\,r^2/2}, \quad r \geq 0.
\en
This classical result can be proved by applying the spherical isoperimetric
inequality (a theorem due to P. L\'evy and Schmidt), or using 
the logarithmic Sobolev inequality on the sphere (due to Mueller and Weissler, cf.
\cite{M-W}, \cite{Le2}). In Section 11 the following assertion is proved.

\vskip5mm
{\bf Theorem 1.3.} {\sl For any function $f$ on $\S^{n-1}$ with
$\|f\|_{\rm Lip} \leq 1$ and $\sigma_{n-1}$-mean $m$,
\be
\sigma_{n-1}\big\{|f-m| \geq r\big\} \leq \frac{1}{r}\,e^{-(n-2)\, r^2/2}\,
\E_{\sigma_{n-1}}\,|f-m|, \quad r \geq 0.
\en
}

\vskip2mm
By the (spherical) Poincar\'e inequality,
the expectation $\E_{\sigma_{n-1}}\,|f-m| = \int_{\S^{n-1}} |f-m|\,d\sigma_{n-1}$
is bounded by $\frac{1}{\sqrt{n-1}}$. Hence (1.8) yields
$$
\sigma_{n-1}\big\{|f-m| \geq r\big\} \leq \frac{1}{r\sqrt{n-1}}\ e^{- (n-2)\,r^2/2}, 
$$
which may be viewed as a sharpening of (1.7), at least for the values
in the range $r \geq 1/\sqrt{n}$. Another point of sharpening is that the factor 
$\E_{\sigma_{n-1}}\,|f-m|$ may be much smaller than $1/\sqrt{n}$.

Covariance representations in Theorems 1.1-1.2 can be further developed in terms 
of derivatives of higher orders such as the spherical Hessian $f''_\S$, which 
is well motivated from the point of view of higher order concentration phenomena.
In particular, it will be shown that on $\S^{n-1} \times \S^{n-1}$ there 
exists a probability measure $\mu_n$ with marginals $\sigma_{n-1}$ such that
\be
{\rm cov}_{\sigma_{n-1}}(f,g) = c_n \int_{\S^{n-1}} \int_{\S^{n-1}} 
\left<Df(x),Dg(y)\right> d\mu_n(x,y)
\en
for all $C^2$-smooth $f,g$ that are orthogonal to all linear functions in 
$L^2(\sigma_{n-1})$, and with constants behaving like $c_n \sim 1/n^2$ for 
the growing dimension $n$. In this representation, we use the linear differential 
operator 
$$
Df(x) = f''_\S(x) - 2\,\nabla_\S f(x) \otimes x
$$ 
and involve 
the canonical inner product for symmetric $n \times n$ matrices
in the integrand in (1.9). Such results are discussed in Sections 12-13.
First applications to bounding the covariance by quantities of order $1/n^2$ 
are given in Section 14, and then a spherical second 
order deviation inequality is derived in Section 15. 

In the end, we give a short overview of one-dimensional covariance
representations and then focus on the case of the circle in Theorem 1.1. 
The whole plan is as follows.

\vskip5mm
1. Introduction

2. Gaussian covariance representation

3. From the Gauss space to the sphere

4. Proof of Theorem 1.1 for $n \geq 3$

5. Asymptotic behaviour of mixing densities

6. Proof of Theorem 1.1 for the circle

7. Differentiation on the sphere

8. Spherical Laplacian

9. Semi-group approach to the covariance

10. Comparison of the two representations

11. Applications to spherical concentration

12. Second order covariance identity in Gauss space

13. Second order covariance identities on the sphere

14. Upper bounds on covariance of order $1/n^2$

15. Second order concentration on the sphere

16. Covariance representations on the line

17. Periodic covariance representations

18. From the circle to the interval

19. Mixing measures on the circle and the interval

\vskip5mm
\section{{\bf Gaussian covariance representation}}
\setcounter{equation}{0}

\vskip2mm
\noindent
The mixing measure in the covariance representation (1.1)
may be constructed as the mixture
\be
\pi_n = \int_0^1 {\sf L}\big(X, tX + \sqrt{1-t^2}\, Z\big)\,dt,
\en
where $X$ and $Z$ are independent standard normal random vectors in $\R^n$ 
with density $\varphi_n$, and where ${\sf L}$ refers to the joint distribution.
Since $tX + \sqrt{1-t^2}\, Z$ and $X$ are distributed according to
$\gamma_n$, it follows from this description that $\pi_n$ has marginals $\gamma_n$.

Let us state the identity (1.1) once more in a way which will be needed later on.

\vskip5mm
{\bf Theorem 2.1.} {\sl Let $u$ and $v$ be $C^1$-smooth complex-valued
functions on an open set $G$ in $\R^n$ of $\gamma_n$-measure $1$ such that
\be
\E\,|\nabla u(X)|^2 < \infty, \quad \E\,|\nabla v(X)|^2 < \infty.
\en
Then the random variables $u(X)$ and $v(X)$ have finite second moments, and
\be
\cov(u(X),v(X)) = \int_{\R^n} \int_{\R^n} 
\left<\nabla u(x),\nabla v(y)\right> d\pi_n(x,y).
\en
}

\vskip2mm
The double integral in (2.3) is bounded in absolute value by
$$
\int_{\R^n} \int_{\R^n} |\left<\nabla u(x),\nabla v(y)\right>|\, d\pi_n(x,y) \leq
\frac{1}{2}\,\E\,|\nabla u(X)|^2 + \frac{1}{2}\,\E\,|\nabla v(X)|^2,
$$
so that it is finite. Also, the Poincar\'e-type inequality
\be
\Var(u(X)) \leq \E\,|\nabla u(X)|^2
\en
(which readily follows from (1.1)) and the assumption (2.2) ensures that 
$u(X)$ and $v(X)$ have finite second moments,
so the covariance is well-defined and finite as well.

The equality (2.3) is easily verified to be true for all exponential functions 
$u(x) = e^{i\left<t,x\right>}$, $v(x) = e^{i\left<s,x\right>}$. Therefore, it 
holds true for finite linear combinations of such functions and actually for 
Fourier-Stieltjes transforms of signed measures on $\R^n$ with finite absolute 
moments. This class is sufficiently rich to properly approximate all smooth 
functions as in Theorem 2.1 and establish (2.3) in full generality 
by an approximation argument.

As we have mentioned before, one may also state the identity (2.3) by
involving the Ornstein-Uhlenbeck operators
\be
P_t v(x) = \int_{\R^n} v\Big(e^{-t} x + \sqrt{1 - e^{-2t}} y\Big)\,d\gamma_n(y),
\quad x \in \R^n, \ t \geq 0.
\en
A similar definition is applied to vector-valued functions.
These operators act on Lebesgue spaces $L^p(\gamma_n)$ as contractions and
form a semi-group, that is, $P_t (P_s v) = P_{t+s} v$, $t,s \geq 0$.

If $v$ is integrable and has an integrable gradient $\nabla v$ with respect to 
$\gamma_n$, one may differentiate under the integral sign in (2.5), which leads to
$$
\nabla P_t v(x) = e^{-t} 
\int_{\R^n} \nabla v\Big(e^{-t} x + \sqrt{1 - e^{-2t}} y\Big)\,d\gamma_n(y) =
e^{-t} P_t \nabla v(x).
$$

An equivalent formulation of (2.3) is the following:

\vskip5mm
{\bf Theorem 2.2.} {\sl Under the assumptions of Theorem $2.1$,
\be
\cov(u(X),v(X)) = \int_0^\infty \E \left<\nabla u(X),\nabla P_t v(X)\right> dt,
\en
where the integral is absolutely convergent.
}

\vskip5mm
From the definition (2.1) it follows that the probability measure $\pi_n$ is 
absolutely continuous with respect to Lebesgue measure on $\R^n \times \R^n$ 
and has some density
$$
p_n(x,y) = \frac{d\pi_n(x,y)}{dx\,dy}.
$$
To write it down explicitly, let us write $s = \sqrt{1-t^2}$. For any bounded 
Borel measurable function $h:\R^n \times \R^n \rightarrow \R$, changing 
the variable $z = \frac{y - tx}{s}$, we have
\bee
\E\, h(X,tX + sZ) 
 & = & 
\int_{\R^n}\!\int_{\R^n}
h(x,tx + sz)\,\varphi_n(x)\varphi_n(z)\ dx\,dz \\
 & = & 
\frac{1}{s^n} \int_{\R^n}\!\int_{\R^n} 
h(x,y)\,\varphi_n(x)\varphi_n\Big(\frac{y-tx}{s}\Big)\ dx\,dy,
\ene
which means that
$$
p_n(x,y) = \int_0^1 s^{-n}
\varphi_n(x)\varphi_n\bigg(\frac{y-tx}{s}\bigg)\,dt.
$$
But
$$
\varphi_n(x)\varphi_n\bigg(\frac{y-tx}{s}\bigg) =
\frac{1}{(2\pi)^n} \ \exp\bigg[-
\frac{|x|^2 + |y|^2 -2t\left<x,y\right>}{2s^2}\bigg].
$$
Thus, we arrive at:

\vskip5mm
{\bf Proposition 2.3.} {\sl The mixing probability measure $\pi_n$ 
in the Gaussian covariance representation $(2.3)$ has density
\be
p_n(x,y) = \frac{1}{(2\pi)^n} \,
\int_0^1 s^{-n} \,\exp\bigg[-
\frac{|x|^2 + |y|^2 -2t\left<x,y\right>}{2s^2}\bigg]\,dt, \quad
x,y \in \R^n,
\en
where we write $s = \sqrt{1-t^2}$.
}

\vskip5mm
If $n \geq 2$, the density $p_n$ is unbounded, since
$p_n(0,0) = (2\pi)^{-n} \, \int_0^1 s^{-n}\,dt = \infty$.
On the other hand, in the case $n=1$, we have $p_1(x,y) \leq \frac{1}{4}$.

\vskip5mm
\section{{\bf From the Gauss space to the sphere}}
\setcounter{equation}{0}

\vskip2mm
\noindent
Starting from (2.3), one can develop a spherical variant of the Gaussian 
covariance identity
\be
{\rm cov}_{\gamma_n}(u,v) = \int_G \int_G
\left<\nabla u(x),\nabla v(y)\right> p_n(x,y)\,dx\,dy.
\en
Here $G$ is an open subset of $\R^n$ of a full $\gamma_n$-measure,
and $p_n$ is the density of the probability measure $\pi_n$  
described in Proposition 2.3. This section is devoted to the first step
in this reduction. Putting
$$
x = r\theta, \ y = r'\theta', \ r,r'>0, \ \theta, \theta' \in \S^{n-1}, 
$$
the formula (2.7) takes the form in polar coordinates
\be
p_n(r\theta,r'\theta') = \frac{1}{(2\pi)^n} \,
\int_0^1 s^{-n} \,\exp\bigg[-
\frac{r^2 + r'^2 -2rr' t\left<\theta,\theta'\right>}{2s^2}\bigg]\,dt,
\en
where we use the notation $s = \sqrt{1-t^2}$ in the integrand, as before.

With two smooth functions $f,g:\S^{n-1} \rightarrow \R$, 
we associate the functions $u(x) = f(\theta)$, $v(y) = g(\theta')$, where
$$
\theta = \frac{x}{r}, \ \theta' = \frac{y}{r'}, \qquad 
r = |x|, \ r' = |y|.
$$
They are defined and $C^1$-smooth on $G = \R^n \setminus \{0\}$ 
and have (Euclidean) gradients
\be
\nabla u(x) = \frac{1}{r}\, \nabla_\S f(\theta), \quad
\nabla v(y) = \frac{1}{r'}\, \nabla_\S g(\theta').
\en

As random variables on the sphere, $f$ and $g$ have the same 
distributions under the normalized Lebesgue measure $\sigma_{n-1}$ 
on $\S^{n-1}$ as $u$ and $v$ have under $\gamma_n$. In particular,
$$
\int_{\R^n}\,u\,d\gamma_n = \int_{\S^{n-1}} f\,d\sigma_{n-1}, \quad
\int_{\R^n}\,v\,d\gamma_n = \int_{\S^{n-1}} g\,d\sigma_{n-1}.
$$
A similar identity is also true for the products $uv$ ang $fg$.
Hence from (3.1) we obtain that
\be
{\rm cov}_{\sigma_{n-1}}(f,g) = {\rm cov}_{\gamma_n}(u,v) =
\int_{\R^n}\!\int_{\R^n} 
\left<\nabla u(x),\nabla v(y)\right> p_n(x,y)\,dx\,dy.
\en

In addition, since $r$ and $\theta$ are independent under $\gamma_n$, and
$r$ has the same distribution as $|Z|$ (where $Z$ is a standard
normal random vector in $\R^n$), we have
$$
\int |\nabla u|^2\,d\gamma_n = \E\,\frac{1}{|Z|^2}\,
\int_{\S^{n-1}} |\nabla_\S f|^2\,d\sigma_{n-1}
$$
and
$$
\int |\nabla v|^2\,d\gamma_n = \E\,\frac{1}{|Z|^2}\,
\int_{\S^{n-1}} |\nabla_\S g|^2\,d\sigma_{n-1}.
$$
Here, while the spherical integrals are finite by the smoothness of $f$
and $g$, the expectation on the right-hand side is finite for 
$n \geq 3$, only. So, this has to be assumed.
For further developments, we will need a recursive formula on
similar expectations. 

\vskip5mm
{\bf Lemma 3.1.} {\sl For any real number $m > -(n-2)$,
$$
\E\,|Z|^m = (m+n-2)\,\E\,|Z|^{m-2},
$$
where $Z$ is a standard normal random vector in $\R^n$.
In particular, $\E\,|Z|^{-2} = \frac{1}{n-2}$ for $n \geq 3$.
}

\vskip5mm
{\bf Proof.} One may use a general formula
$$
\E\,w(|Z|) \, = \,
\frac{1}{(2\pi)^{n/2}} \int_{\R^n} w(|x|)\,e^{-|x|^2/2}\,dx \, = \,
\frac{\omega_{n-1}}{(2\pi)^{n/2}} \int_0^\infty w(r) r^{n-1}\,e^{-r^2/2}\,dr,
$$
where we integrated in polar coordinates in the last step. Here and in the sequel, 
we denote by $\omega_{n-1}$ the surface area of the sphere $\S^{n-1}$, that is,
\be
\omega_{n-1} = 
\frac{(2\pi)^{\frac{n}{2}}}{2^{\frac{n}{2} - 1}\, \Gamma(\frac{n}{2})}.
\en
Integrating by parts for the polynomial function $w(r) = r^m$ with
$m > - (n-2)$ (this condition is posed for the integrability reason), we get
\bee
\E\,|Z|^m 
 & = &
\frac{\omega_{n-1}}{(2\pi)^{n/2}} \int_0^\infty r^{m+n-1}\,e^{-r^2/2}\,dr \\
 & = &
(m+n-2)\,\frac{\omega_{n-1}}{(2\pi)^{n/2}} 
\int_0^\infty r^{m+n-3}\,e^{-r^2/2}\,dr \, = \, (m+n-2)\,\E\,|Z|^{m-2}.
\ene
\qed

\vskip5mm
\section{{\bf Proof of Theorem 1.1 for $n \geq 3$}}
\setcounter{equation}{0}

\vskip2mm
\noindent
Keeping the same notations as before, let us return to (3.4).
To continue, one may integrate on the right-hand of this representation 
in polar coordinates using the general formula 
\bee
\int_{\R^n}\!\int_{\R^n} h(x,y)\,dx dy
 & = & 
\int_0^\infty\!\int_0^\infty \int_{\S^{n-1}}\int_{\S^{n-1}}
h(r\theta,r'\theta')\,(r r')^{n-1} \,dr\, dr' \ d\theta\, d\theta' \\
 & \hskip-20mm = &
 \hskip-10mm 
\omega_{n-1}^2 \int_0^\infty\!\int_0^\infty \int_{\S^{n-1}}\int_{\S^{n-1}}
h(r\theta,r'\theta')\,(r r')^{n-1} \,dr\, dr' \
d\sigma_{n-1}(\theta)\, d\sigma_{n-1}(\theta'), 
\ene
where $\omega_{n-1}$ stands for the area size of the sphere $\S^{n-1}$
described in (3.5). Hence, applying (3.2)-(3.3), the double integral 
in (3.4) may be rewritten, and we obtain that
\be
{\rm cov}_{\sigma_{n-1}}(f,g) = \int_{S^{n-1}}\int_{\S^{n-1}} 
\left<\nabla_\S f(\theta),\nabla_\S g(\theta')\right>
\psi_n(\theta,\theta')\,d\sigma_{n-1}(\theta)\, d\sigma_{n-1}(\theta')
\en
with
\bee
\psi_n(\theta,\theta') 
 & = &
\frac{1}{2^{n - 2}\, \Gamma(\frac{n}{2})^2}
\int_0^1 s^{-n} \bigg[\int_0^\infty\!\int_0^\infty 
\exp\bigg[-\frac{r^2 + r'^2 -2rr' t\left<\theta,\theta'\right>}{2s^2}\bigg]
(r r')^{n-2} \,dr\, dr'\bigg] dt \\
 & = &
\frac{1}{2^{n - 2}\, \Gamma(\frac{n}{2})^2}
\int_0^1 s^{n-2} \, \bigg[\int_0^\infty\!\int_0^\infty 
\exp\bigg[-\frac{r^2 + r'^2 -2rr' t\left<\theta,\theta'\right>}{2}\bigg]
(r r')^{n-2} \,dr\, dr'\bigg] dt,
\ene
where changed the variables in the last step by replacing $r$ with $sr$ and
$r'$ with $s r'$.

Let us see that $\psi_n$ is integrable over the product measure
$\sigma_{n-1} \otimes \sigma_{n-1}$, so that it serves as density
of some finite Borel positive measure $\nu_n$ on $\S^{n-1} \times \S^{n-1}$
with total mass
$$
c_n = \nu_n(\S^{n-1} \times \S^{n-1}) = 
\int_{\S^{n-1}}\int_{\S^{n-1}} \psi_n(\theta,\theta')\,
d\sigma_{n-1}(\theta)\,d\sigma_{n-1}(\theta').
$$
Repeating integration in polar coordinates as before, 
one may notice that
$$
\E \int_0^1 \frac{1}{|X|\,|tX + sZ|}\,dt
=
\int_{\R^n}\!\int_{\R^n} \frac{1}{|x|\,|y|}\, p_n(x,y)\,dx\, dy = c_n.
$$
But, for any fixed $0 < t < 1$, by Cauchy's inequality and applying Lemma 3.1,
\bee
\E\, \frac{1}{|X|\,|tX + sZ|} 
 & < &
\bigg(\E\, \frac{1}{|X|^2}\bigg)^{1/2}\,
\bigg(\E\, \frac{1}{|tX + sZ|^2}\bigg)^{1/2} \\
 & = &
\E\, \frac{1}{|X|^2} \, = \, \frac{1}{n-2} \, < \, \infty.
\ene
Thus, $\nu_n = c_n \mu_n$ for some Borel probability measure $\mu_n$ on 
$\S^{n-1} \times \S^{n-1}$ with $c_n < \frac{1}{n-2}$.

Note also that $\psi_n(\theta,\theta') = \psi_n(\left<\theta,\theta'\right>)$ 
depends on $\theta$ and $\theta'$ via the inner product only.
This implies that the integral
$$
\int_{\S^{n-1}} \psi_n(\theta,\theta')\,d\sigma_{n-1}(\theta') =
\int_{\S^{n-1}} \psi_n(\left<\theta,\theta'\right>)\,
d\sigma_{n-1}(\theta')
$$
does not depend on $\theta$ and is therefore equal to $c_n$. 
As a consequence, we conclude that $\nu_n$ is a positive, finite, symmetric 
measure whose marginal distributions are spherically invariant measures.
Hence they coincide with a multiple of $\sigma_{n-1}$. Returning to (4.1),
one can summarize
in the following refinement of Theorem 1.1 for the case $n \geq 3$.

\vskip5mm
{\bf Theorem 4.1.} {\sl For all smooth functions $f,g$ on $\S^{n-1}$
$(n \geq 3)$,
\be
{\rm cov}_{\sigma_{n-1}}(f,g) = \int_{\S^{n-1}}\int_{\S^{n-1}} 
\left<\nabla_\S f(x),\nabla_\S g(y)\right> d\nu_n(x,y),
\en
where $\nu_n$ is a positive measure on $\S^{n-1} \times \S^{n-1}$ with
density $\psi_n(\left<x,y\right>)$ with respect to the product measure 
$\sigma_{n-1} \otimes \sigma_{n-1}$, where
\be
\psi_n(\alpha) =
\frac{1}{2^{n - 2}\, \Gamma(\frac{n}{2})^2}
\int_0^1 s^{n-2} \, \bigg[\int_0^\infty\!\int_0^\infty 
\exp\bigg[-\frac{r^2 + r'^2 -2rr' t\alpha}{2}\bigg]
(r r')^{n-2} \,dr\, dr'\bigg] dt
\en
defined for $|\alpha| \leq 1$. Moreover, $\nu_n = c_n \mu_n$ for 
a probability measure $\mu_n$ on $\S^{n-1} \times \S^{n-1}$ with marginals 
$\sigma_{n-1}$ and with total mass satisfying
\be
\frac{1}{n-1} < c_n < \frac{1}{n-2}.
\en
}

\vskip2mm
{\bf Proof.} It remains to explain the lower bound in (4.3).
From (4.2) in the case of equal functions $f = g$ we get
\be
\Var(f(\xi)) = c_n \int_{\S^{n-1}}\int_{\S^{n-1}} 
\left<\nabla_\S f(\theta),\nabla_\S f(\theta')\right>
d\mu_n(\theta,\theta'),
\en
where $\xi$ is a random vector with distribution $\sigma_{n-1}$. By Cauchy's inequality,
\begin{eqnarray}
\Var(f(\xi))
 & \leq &
c_n \bigg(\int\!\!\int
|\nabla_\S f(x)|^2\,d\mu_n(x,y)\bigg)^{1/2} 
\bigg(\int\!\!\int |\nabla_\S f(y)|^2\,d\mu_n(x,y)\bigg)^{1/2}\\
 & = &
c_n \int |\nabla_\S f|^2\,d\sigma_{n-1}. \nonumber
\end{eqnarray}
Thus, we arrive at the Poincar\'e-type inequality 
$\Var(f(\xi)) \leq c_n\, \E\, |\nabla_\S f(\xi)|^2$ on the unit sphere.
Here, the constant may not be better than the optimal constant $\frac{1}{n-1}$.

Moreover, the equality in (4.6) is only possible if and only if
$$
\left<\nabla_\S f(x),\nabla_\S f(y)\right> = |\nabla_\S f(x)|\,|\nabla_\S f(y)|
$$
for $\mu_n$-almost all $(x,y) \in \S^{n-1} \times \S^{n-1}$, hence for all 
$x,y \in \S^{n-1}$ (by the continuity of the gradient). Equivalently, 
one should have $\nabla_\S f(x) = \nabla_\S f(y)$ for all $x,y \in \S^{n-1}$,
which means that $f$ must be a constant. This implies that $c_n > \frac{1}{n-1}$.
\qed

\vskip5mm
{\bf Remark 4.2.} In order to estimate $c_n$ and prove (4.4), it is also 
sufficient to apply the identity (4.5) to linear functions 
$f(x) = \left<x,v\right>$, $v \in \S^{n-1}$. In this case
$\Var(f(\xi)) = \frac{1}{n}$ and
$\nabla_\S f(\theta) = v - \left<v,\theta\right>\theta$, so that
$$
\left<\nabla_\S f(\theta),\nabla_\S f(\theta')\right> =
1 - \left<v,\theta\right>^2  - \left<v,\theta'\right>^2 +
\left<v,\theta\right>\left<v,\theta'\right> \left<\theta,\theta'\right>.
$$
Hence, integrating this expression over $\mu_n$ and using the property
that $\mu_n$ has marginals $\sigma_{n-1}$, (4.5) yields
$$
\frac{1}{nc_n} = 1 - \frac{2}{n} +
\int_{\S^{n-1}}\int_{\S^{n-1}} 
\left<v,\theta\right>\left<v,\theta'\right> \left<\theta,\theta'\right>
d\mu_n(\theta,\theta').
$$
Here, the left-hand side does not depend on $v$, while formally
the right-hand side does. Integrating this equality over 
$d\sigma_{n-1}(v)$, we then get
$$
\frac{1}{nc_n} = 1 - \frac{2}{n} + \frac{1}{n}\,
\int_{\S^{n-1}}\int_{\S^{n-1}} 
\left<\theta,\theta'\right>^2 d\mu_n(\theta,\theta').
$$
The latter implies $\frac{1}{nc_n} > 1 - \frac{2}{n}$, that is, 
$c_n < \frac{1}{n-2}$. On the other hand, since 
$\left<\theta,\theta'\right>^2 < 1$ for $\mu_n$-almost all couples
$(\theta,\theta')$, the double integral is smaller than 1. Hence,
$\frac{1}{nc_n} < 1 - \frac{1}{n}$, that is, we also have a lower bound
$c_n > \frac{1}{n-1}$.


\vskip5mm
\section{{\bf Asymptotic behaviour of mixing densities}}
\setcounter{equation}{0}

\vskip2mm
\noindent
Since $\left<\theta,\theta'\right> = 1 - |\theta - \theta'|^2$,
the property that the density function $\psi_n(\theta,\theta')$ on 
$\S^{n-1} \times \S^{n-1}$ depends on the inner product is equivalent to 
the assertion that it depends on the distance function $|\theta - \theta'|$. 
We will now examine an asymptotic behaviour of this function when
this distance is close to zero.

The previous results are applicable in dimension $n=2$ as well 
(with a slightly modified argument), that is, for the circle $\S^1$ 
which we equip with the uniform probability measure $\sigma_1$. In this case, 
according to (4.3), the measure $\nu_2$ has density 
\be
\psi_2(\theta,\theta') = \psi_2(\left<\theta,\theta'\right>) = 
\int_0^1\, \bigg[\int_0^\infty\!\int_0^\infty 
\exp\bigg[-\frac{r^2 + r'^2 -2rr' t\left<\theta,\theta'\right>}{2}\bigg]\,
dr dr'\bigg] \,dt,
\en
where we use the notation $s = \sqrt{1-t^2}$. In the previous section,
it was however not checked, whether $\nu_2$ is finite or not.
In fact, although the above triple integral is apparently not easy to
evaluate explicitly, one can show that the density $\psi_2$ is bounded.

Since $\sigma_1 \otimes \sigma_1$ is invariant under
transformations $(\theta,\theta') \rightarrow (\pm \theta,\pm\theta')$, we have
$$
\nu_2(\S^1 \times \S^1) = \frac{1}{4} \int_{\S^1}\int_{\S^1} 
\Psi(\alpha)\,d\sigma_1(\theta)\,d\sigma_1(\theta'), \quad 
\alpha = \left<\theta,\theta'\right>,
$$
where
$$
\Psi(\alpha) = \int_0^1 
\bigg[\int_{-\infty}^\infty \int_{-\infty}^\infty 
\exp\bigg[-\frac{r^2 + r'^2 -2rr' t\alpha}{2}\bigg]\,dr dr'\bigg] \,dt.
$$
In addition,
\be
\psi_2(\alpha) \leq \Psi(\alpha).
\en

The inner integral with respect to $r$ in the definition of $\Psi(\alpha)$ 
is equal to 
$$
\int_{-\infty}^\infty 
\exp\bigg[-\frac{(r -r' t\alpha)^2 - (r' t\alpha)^2 + r'^2}{2}\bigg]\,dr =
\sqrt{2\pi}\, \exp\bigg[-\frac{1 - (t\alpha)^2}{2}\,r'^2\bigg].
$$
The next integration with respect to $r'$ yields the value
$\frac{2\pi}{\sqrt{1 - (t\alpha)^2}}$.
Integrating over $0<t<1$, we then get
$$
\Psi(\alpha) = 2\pi \int_0^1 \frac{dt}{\sqrt{1 - (t\alpha)^2}} =
2\pi\, \frac{\arcsin(\alpha)}{\alpha}.
$$
In view of the concavity of the sine function on $[0,\frac{\pi}{2}]$,
we have $\frac{\pi}{2}\,\sin(t) \geq t$ on this interval, or equivalently
$\arcsin(\alpha) \leq \frac{\pi}{2}\, \alpha$ for all $0 \leq \alpha \leq 1$.
Thus, $\Psi(\alpha) \leq \pi^2$ implying the following:

\vskip5mm
{\bf Lemma 5.1.} {\sl We have $\psi_2(\theta,\theta') \leq \pi^2$ for all
$\theta,\theta' \in \S^1$ and 
$$
c_2 = \nu_2(\S^1 \times \S^1) \leq \frac{\pi^2}{4}.
$$
}

However, the density of the mixing measure $\nu_n$ is not bounded for 
$n \geq 3$. Let us give natural lower and upper bounds for 
$\psi_n(\alpha)$ in terms of $\alpha$ and $n$ as $\alpha$ approaches 1. 
To this aim, consider the functions
\be
R_m(\alpha) = \int_0^1 s^m \bigg[\int_0^\infty \!\!\int_0^\infty 
\exp\bigg[-\frac{r^2 + r'^2 -2rr' t\alpha}{2}\bigg]\,(r r')^m\,dr\, dr'\bigg] \,dt,
\quad 
\en
for $-1 \leq \alpha < 1$
with a real parameter $m \geq 0$ which may be used with $m = n-2$
for the sphere $\S^{n-1}$, $n \geq 3$, and with $0<m<1$ for $n=2$. 

Put $\beta = 1-\alpha t$ and first consider the inner double integral
\bee
I_m(\beta) 
 & = &
\int_0^\infty \!\!\int_0^\infty 
\exp\bigg[-\frac{r^2 + r'^2 -2(1-\beta) rr'}{2}\bigg]\,(r r')^m\,dr\, dr' \\
 & = &
\int_0^\infty \!\!\int_0^\infty 
\exp\bigg[-\frac{1}{2}\,(r - r')^2 - \beta r r'\bigg]\,(r r')^m\,dr\, dr'.
\ene
Changing the variables $x = r - r'$, $y = r r'$, the last 
integral takes the form
$$
\int_0^\infty \!\!\int_0^\infty 
\exp\bigg[-\frac{1}{2}\,x^2 - \beta y\bigg]\,y^m\,dr dr'.
$$
The next change $r' = y/x$, $x = r - y/x$ with arbitrary $x \in \R$ and $y>0$
leads to
\be
I_m(\beta) = \int_{-\infty}^\infty \int_0^\infty 
\exp\bigg[-\frac{1}{2}\,x^2 - \beta y\bigg]\,
\frac{y^m}{\sqrt{x^2 + 4y}}\,dx\, dy,
\en
which is bounded from above by
$$
\frac{1}{2}
\int_{-\infty}^\infty \int_0^\infty \exp\bigg[-\frac{1}{2}\,x^2 - \beta y\bigg]\,
y^{m - 1/2}\,dx\, dy \, = \,
\frac{\sqrt{2\pi}}{2}\,\beta^{-m - 1/2}\ \Gamma(m + 1/2).
$$

For a similar lower bound, one may restrict integration in (5.4)
to the regions $|x| < 1$ and $y > 1$, in which case $x^2 + 4y \leq 5y$.
Since also $0 < \beta < 1$, this leads to
\bee
I_m(\beta) 
 & \geq  &
\frac{2}{\sqrt{5e}} \int_1^\infty e^{- \beta y}\,y^{m-1/2}\,dy = 
\frac{2}{\sqrt{5e}} \,\beta^{-m-1/2} 
\int_\beta^\infty e^{-z}\,z^{m-1/2}\,dz \\
 & \geq &
b_0 \,\beta^{-m - 1/2}\ \Gamma(m + 1/2)
\ene
for some absolute constant $b_0>0$. Thus,
\be
b_0 \,\beta^{-m - 1/2}\ \Gamma(m + 1/2) \leq I_m(\beta) \leq
b_1 \,\beta^{-m - 1/2}\ \Gamma(m + 1/2),
\en
where one may take $b_1 = \frac{\sqrt{2\pi}}{2}$ and
where we recall that $\beta = 1-\alpha t$.

Let us turn to (5.3) and apply the upper bound in (5.5) which gives
\begin{eqnarray}
R_m(\alpha) 
 & = &
\int_0^1 s^m\,I_m(1-\alpha t)\,dt \nonumber \\
 & \leq &
\frac{\sqrt{2\pi}}{2}\ \Gamma(m + 1/2)
\int_0^1 s^m\,(1-\alpha t)^{-m - 1/2}\,dt.
\end{eqnarray}
Note that $R_m(\alpha)$ is increasing in $-1 \leq \alpha \leq 1$.
Using $s \leq \sqrt{2}\,(1 - t)^{1/2}$ 
and assuming $0 \leq m < 1$, this gives
$$
R_m(\alpha) \leq R_m(1) \leq \sqrt{\pi}\ \Gamma(m + 1/2)
\int_0^1 (1 - t)^{-\frac{m + 1}{2}}\,dt = 
\sqrt{\pi}\ \Gamma(m + 1/2)\,\frac{2}{1-m}.
$$
That is, we proved:

\vskip5mm
{\bf Lemma 5.2.} {\sl For $0 \leq m < 1$ and $|\alpha| \leq 1$,
\be
R_m(\alpha) \leq 
2\sqrt{\pi}\ \frac{\Gamma(m + 1/2)}{1-m}.
\en
In particular, choosing $m=0$, we have  
$\psi_2(\theta,\theta') = R_0(\alpha) \leq 2\pi$ for all 
$\theta,\theta' \in \S^1$.
}

\vskip5mm
The last bound is slightly better then the uniform bound of Lemma 5.1.

However, Lemma 5.2 is not applicable for $m \geq 1$, and
we have to return to the bound (5.6). We need to estimate
the integral
$$
J_m(\alpha) = \int_0^1 s^m \,(1 - \alpha t)^{-m - 1/2}\,dt
$$
as $\alpha$ approaching 1. Let us write $Q_1 \sim Q_2$
for two positive quantities, if the ratio $Q_1/Q_2$ is bounded
away from zero and infinity by some $m$-dependent constants.
Note that 
$J_m(\alpha) \leq J_m(1/2) \leq 2^{m+1/2}$ for $\alpha \leq 1/2$,
so, we may assume that $\frac{1}{2} \leq \alpha < 1$.
Since the above integral remains bounded when integrating over 
$0<t<1/2$, and $s \sim (1 - t)^{1/2}$, we have, changing 
the variable and putting $\ep = 1-\alpha$,
\bee
J_m(\alpha) 
 & \sim &
\int_0^{1/2} \frac{t^{m/2}}{(1 - \alpha(1-t))^{m + 1/2}}\,dt = 
\int_0^{1/2} \frac{t^{m/2}}{(\ep + (1-\ep)\, t)^{m + 1/2}}\,dt \\
 & = &
\ep^{-\frac{m-1}{2}}
\int_0^{1/2\ep} \frac{x^{m/2}}{(1 + (1-\ep)\,x)^{m + 1/2}}\,dt
 \sim
\ep^{-\frac{m-1}{2}}
\int_0^{1/2\ep} \frac{x^{m/2}}{(1 + x)^{m + 1/2}}\,dx \\
 & \sim &
\ep^{-\frac{m-1}{2}}
\int_{1/2}^{1/2\ep} \frac{1}{(1 + x)^{\frac{m + 1}{2}}}\,dx.
\ene
The last integral is bounded for $m>1$ and is equivalent to
$\log(1/\ep)$ for $m=1$. Therefore, 
$$
R_m(\alpha) \sim (1-\alpha)^{-\frac{m-1}{2}} \ \ (m>1), \qquad
R_1(\alpha) \sim \log\frac{1}{1 - \alpha}.
$$
Using this with $m=n-2$ according to the formula (4.3), 
we obtain a similar conclusion about the densities 
$\psi_n(\alpha) = \psi_n(\left<\theta,\theta'\right>)$ of $\nu_n$
and thus prove the asymptotic relations (1.5):

\vskip5mm
{\bf Proposition 5.3.} {\sl For $-1 \leq \alpha < 1$,
$$
\psi_n(\alpha) \sim (1-\alpha)^{-\frac{n-3}{2}} \ \ (n \geq 4), \qquad
\psi_3(\alpha) \sim \log\frac{1}{1 - \alpha}.
$$
}

\vskip5mm
\section{{\bf Proof of Theorem 1.1 for the circle}}
\setcounter{equation}{0}

\vskip2mm
\noindent
Here we extend the assertion of Theorem 4.1 to the case of the circle,
using the same formula (4.3) for the density $\psi_2$ for the mixing measure
$\nu_2$ which is also described in (5.1).

For this aim, the Gaussian covariance identity may still be used, however, 
with slightly different functions $u$ and $v$.
Namely, with a given smooth function $f:\S^{n-1} \rightarrow \R$, 
we associate the homogeneous functions of order $\ep > 0$
$$
u_\ep(x) = r^\ep f(r^{-1} x) = r^\ep f(\theta), \quad 
{\rm where} \quad
r = |x| = \sqrt{x_1^2 + \dots + x_n^2}, \ \ \theta = r^{-1} x.
$$
They are defined, $C^1$-smooth in the whole space 
except for the origin, and have gradients
$$
\nabla u_\ep(x) = r^{\ep - 1}\, 
\big(\ep f(\theta) \theta + \nabla_\S f(\theta)\big).
$$
Define similarly the homogeneous functions $v_\ep$ for 
a smooth function $g:\S^{n-1} \rightarrow \R$, so that
$$
\nabla v_\ep(x) = r^{\ep - 1}\, 
\big(\ep g(\theta) \theta + \nabla_\S g(\theta)\big).
$$
Since $r$ and $\theta$ are independent under $\gamma_n$, 
while $r$ has the same distribution as $|Z|$, we have
\be
{\rm cov}_{\gamma_n}(u_\ep,v_\ep) =
\E\,|Z|^{2\ep}\, {\rm cov}_{\sigma_{n-1}}(f,g) \, \rightarrow \,
{\rm cov}_{\sigma_{n-1}}(f,g)
\en
as $\ep \rightarrow 0$. In addition,
$$
\int |\nabla u_\ep|^2\,d\gamma_n = \E\,|Z|^{2\ep - 2}\,
\int_{\S^{n-1}} (\ep^2 f^2 + |\nabla_\S f|^2)\,d\sigma_{n-1},
$$
which is finite, and similarly for $v_\ep$. Hence, we may apply 
Theorem 2.1 to the couple $(u_\ep,v_\ep)$, which gives
\be
{\rm cov}_{\gamma_n}(u_\ep,v_\ep) = \int_{\R^n}\int_{\R^n} 
\left<\nabla u_\ep(x),\nabla v_\ep(y)\right> p_n(x,y)\,dx\,dy.
\en

Putting
\bee
w_\ep(\theta,\theta') 
 & = &
\ep^2 f(\theta)g(\theta') \left<\theta,\theta'\right> +
\ep f(\theta) \left<\theta,\nabla_\S g(\theta')\right> \\
 & & + \  
\ep g(\theta') \left<\theta',\nabla_\S f(\theta)\right> +
\left<\nabla_\S f(\theta),\nabla_\S g(\theta')\right>,
\ene
one may integrate in polar coordinates separately with respect to $x$ and $y$
in (6.2) and represent the double integral 
similarly as in the proof of Theorem 4.1 as
\be
\int_{\S^{n-1}}\int_{\S^{n-1}} w_\ep(\theta,\theta') \,
\psi_{n,\ep}(\theta,\theta')\,d\sigma_{n-1}(\theta) d\sigma_{n-1}(\theta')
\en
with $\psi_{n,\ep}(\theta,\theta')$ given by
$$
\frac{1}{2^{n - 2}\, \Gamma(\frac{n}{2})^2}
\int_0^1 s^{-n} \bigg[\int_0^\infty\!\int_0^\infty 
\exp\bigg[-\frac{r^2 + r'^2 -2rr' t\left<\theta,\theta'\right>}{2s^2}\bigg]
(r r')^{n + \ep - 2} \,dr dr'\bigg] dt,
$$
where we use the notation $s = \sqrt{1-t^2}$. Changing the variable and using 
$\alpha = \left<\theta,\theta'\right>$, this expression is simplified to
$$
\frac{1}{2^{n - 2}\, \Gamma(\frac{n}{2})^2}
\int_0^1 s^{n + 2\ep - 2} \, \bigg[\int_0^\infty\!\int_0^\infty 
\exp\bigg[-\frac{r^2 + r'^2 -2rr' t\alpha}{2}\bigg]
(r r')^{n+\ep-2} \,dr dr'\bigg] dt.
$$

As we know from Section 4, one may continue with $\ep = 0$
in the case $n \geq 3$. If $n=2$, then
$$
\psi_{2,\ep}(\alpha) = \int_0^1 s^{2\ep} \, \bigg[\int_0^\infty\!\int_0^\infty 
\exp\bigg[-\frac{r^2 + r'^2 -2rr' t\alpha}{2}\bigg]
(r r')^\ep \,dr dr'\bigg] dt,
$$
and by (6.2),
\be
{\rm cov}_{\gamma_2}(u_\ep,v_\ep) = 
\int_{\S^1}\int_{\S^1} w_\ep(\theta,\theta') \,
\psi_{2,\ep}(\theta,\theta')\,d\sigma_1(\theta) d\sigma_1(\theta').
\en
Here, the covariance part is convergent as $\ep \rightarrow 0$ to 
${\rm cov}_{\sigma_1}(f,g)$, according to (6.1).

In order to turn to the limit on the right-hand side, note that the 
functions $w_\ep(\theta,\theta')$ remain bounded uniformly over all 
$\theta,\theta' \in \S^1$, $0 < \ep \leq 1/4$, and are convergent to 
$\left<\nabla_\S f(\theta),\nabla_\S g(\theta')\right>$ as $\ep \rightarrow 0$.
By Lemma 6.2 with $m = \ep$, the quantities
$\psi_{2,\ep}(\alpha) = R_{2\ep}(\alpha)$ are also bounded by a constant.
Hence, if we show that $\psi_{2,\ep}(\alpha) \rightarrow \psi_2(\alpha)$
for all $\alpha \in [-1,1)$, one may apply the Lebesgue dominated convergence
theorem, so that to derive from (6.4) the limit case
$$
{\rm cov}_{\sigma_1}(f,g) = 
\int_{\S^1}\int_{\S^1} \left<\nabla_\S f(\theta),\nabla_\S g(\theta')\right>
\psi_2(\theta,\theta')\,d\sigma_1(\theta) d\sigma_1(\theta').
$$

As we did in the previous section, cf. (5.3)-(5.4) with $m=2\ep$, we have
$$
\psi_{2,\ep}(\alpha) = R_{2\ep}(\alpha) =
\int_0^1
\int_{-\infty}^\infty \int_0^\infty 
\exp\bigg[-\frac{1}{2}\,x^2 - (1 - t\alpha) y\bigg]\,
\frac{(s y)^{2\ep}}{\sqrt{x^2 + 4y}}\,dx\, dy\,dt,
$$
of which the limit case is given by
\be
\psi_2(\alpha) = R_0(\alpha) =
\int_0^1
\int_{-\infty}^\infty \int_0^\infty 
\exp\bigg[-\frac{1}{2}\,x^2 - (1 - t\alpha) y\bigg]\,
\frac{1}{\sqrt{x^2 + 4y}}\,dx\, dy\,dt.
\en
In order to apply the Lebesgue dominated convergence
theorem and show that $\psi_{2,\ep}(\alpha) \rightarrow \psi_2(\alpha)$,
it is sufficient to see that the integrands
$$
K_{\ep,\alpha}(x,y,t) = \exp\bigg[-\frac{1}{2}\,x^2 - (1 - t\alpha) y\bigg]\,
\frac{(sy)^{2\ep}}{\sqrt{x^2 + 4y}}
$$
are bounded by an integrable function $K_\alpha(x,y,t)$ on
$\R \times (0,\infty) \times (0,1)$ with respect to the
Lebesgue measure on this product space uniformly over all
$0 < \ep \leq \frac{1}{4}$. Clearly, for a fixed value of
$\alpha \in [-1,1)$, one should examine the majorant
$$
K_\alpha(x,y,t) \, = \sup_{0 < \ep \leq \frac{1}{4}} 
K_{\ep,\alpha}(x,y,t) \, = \,
\exp\bigg[-\frac{1}{2}\,x^2 - (1 - t\alpha) y\bigg]\,
\frac{\max(1,\sqrt{sy})}{\sqrt{x^2 + 4y}}.
$$
Clearly, the last expression is maximized for $\alpha = 1$.
Since also $s \leq 2(1-t)$, we have
$$
K_\alpha(x,y,t) \leq K_1(x,y,t) \leq 
\exp\bigg[-\frac{1}{2}\,x^2 - (1 - t) y\bigg]\,
\frac{\max(1,\sqrt{(1-t) y})}{\sqrt{2y}},
$$
so,
\bee
\int_0^1 \int_{-\infty}^\infty \int_0^\infty 
K_\alpha(x,y,t)\,dx\, dy\,dt
 & \leq &
\sqrt{2\pi} \int_0^1 \int_0^\infty 
e^{-(1 - t) y}\,\frac{\max(1,\sqrt{(1-t)y})}{\sqrt{2y}}\ dy\,dt \\
 & = &
\sqrt{2\pi} \int_0^1 \int_0^\infty 
\frac{1}{\sqrt{1-t}} \ 
e^{-z}\,\frac{\max(1,\sqrt{z})}{\sqrt{2z}}\ dz\,dt \ < \ \infty.
\ene
One may summarize.

\vskip5mm
{\bf Proposition 6.1.} {\sl On the torus $\S^1 \times \S^1$ there 
exists a positive measure $\nu_2$ with marginals $c\sigma_1$
for some constant $c>0$ such that, for all smooth functions 
$f,g$ on $\S^1$, we have the covariance representation
\be
{\rm cov}_{\sigma_1}(f,g) = \int_{\S^1}\int_{\S^1} 
\left<\nabla_\S f(\theta),\nabla_\S g(\theta')\right>
d\nu_2(\theta,\theta').
\en
Moreover, the mixing measure $\nu_2$ has density 
$\psi_2(\left<\theta,\theta'\right>)$ for the function $\psi_2$ with respect to 
the product measure $\sigma_1 \otimes \sigma_1$ given in $(6.5)$. We also have 
$\psi_2 \leq 2\pi$ and $c = \nu_2(\S^1 \times \S^1) < \pi$.
}

\vskip5mm
The triple integral in (6.5) may be simplified by changing the variable
$y = \frac{1}{4}\,x^2 z$ and then integrating over $x$ and $t$. 
Then it becomes the one-dimensional integral
$$
\psi_2(\alpha) = \frac{2}{\alpha} \int_0^\infty 
\Big(\log(1 + z/2) - \log(1 + (1-\alpha)z/2)\Big)\,\frac{dz}{z\sqrt{1 + z}}.
$$

\vskip5mm
\section{{\bf Differentiation on the sphere}}
\setcounter{equation}{0}

\vskip2mm
\noindent
Before turning to Theorem 1.2 and applications of the spherical covariance 
representation,  let us recall a few basic formulas related to
differentiation on the sphere.

The spherical gradient $w = \nabla_\S f(\theta)$ of
a smooth function $f$ on $\S^{n-1}$ may be defined as the shortest vector in $\R^n$ 
satisfying the Taylor formula (1.3). If $f$ is defined and smooth on 
the whole space $\R^n$ or a smaller open subset $G$ containing the unit 
sphere, its spherical gradient at a point $\theta \in \S^{n-1}$
is related to the usual gradient $\nabla f(\theta)$
by
\be
\nabla_\S f(\theta) = \nabla f(\theta) - 
\left<\nabla f(\theta),\theta\right>\theta = P_{\theta^\bot} \nabla f(\theta), 
\en
where $P_{\theta^\bot}$ denotes the projection operator in $\R^n$ onto the orthogonal complement of the line containing $\theta$ (the shifted tangent space).
Alternatively, one may start with a smooth function $f:\S^{n-1} \rightarrow \R^n$ 
and extend it to $G = \R^n \setminus \{0\}$, for example, by the formula
\be
u(x) = f(r^{-1} x), \quad r = |x|.
\en
Then we obtain a smooth function on $G$ with
\be
\nabla u(x) = r^{-1}\,\nabla_\S f(r^{-1} x),
\en
which coincides with the spherical gradient when it is restricted to the sphere.

From (7.1) it follows that the spherical derivative $\nabla_\S$ is a linear
operator satisfying the usual rule of differentiation of products
\be
\nabla_\S (fg) = f \nabla_\S g + g \nabla_\S f.
\en

Here and in the sequel, we endow the real linear space of $n \times n$ matrices
with the standard inner product 
$$
\left<A,B\right> = \sum_{i,j=1}^n a_{ij} b_{ij} = {\rm Tr}(AB), \quad
A = (a_{ij})_{i,j=1}^n, \ B = (a_{ij})_{i,j=1}^n,
$$
(where ${\rm Tr}$ means the trace). The associated 
Hilbert-Schmidt norm is
$$
\|A\|_{\rm HS} = \sqrt{\left<A,A\right>} = \Big(\sum_{i,j=1}^n a_{ij}^2\Big)^{1/2}.
$$

Given a $C^2$-smooth function $f$ on $\S^{n-1}$, its second derivative
or the Hessian $B = f''_S(\theta)$ at the point $\theta \in \S^{n-1}$ may
be defined intrinsically using the Taylor formula 
\begin{eqnarray}
f(\theta') - f(\theta) 
 & = & 
\left<\nabla_\S f(\theta),\theta' - \theta\right> \nonumber \\ 
 & & \hskip-15mm + \
\frac{1}{2}\,\left<B (\theta' - \theta),\theta' - \theta\right> +
o(|\theta' - \theta|^2), \quad \theta' \rightarrow \theta \ \ 
(\theta' \in \S^{n-1}),
\end{eqnarray} 
where $B$ is required to be a symmetric $n \times n$ matrix with 
the smallest Hilbert-Schmidt norm. One can show that, for all 
$v \in \R^n$ and $\theta \in \S^{n-1}$,
\be
f''_\S(\theta) v = \nabla_\S \left<\nabla_\S f(\theta),v\right> + 
\left<v,\theta\right> \nabla_\S f(\theta),
\en
which means that the second order derivative represents the double 
application of the spherical differentiation, however, for $v$ in 
the tangent space, only. 

For a matrix description in terms of (any) smooth 
extension $u$ of $f$ from the sphere to its neighborhood, we have 
\be
f_\S''(\theta) = P_{\theta^\bot} A P_{\theta^\bot}, \quad
A = u''(\theta) - \left<\nabla u(\theta),\theta\right> I_n,
\en
where $I_n$ denotes the identity $n \times n$ matrix, and where $u''$ is 
the usual (Euclidean) Hessian of $u$,  that is, the matrix of second 
order partial derivatives of $u$ (cf. e.g. \cite{B-C-G2}). 

In view of the orthogonality of $\nabla_\S f(\theta)$ and $\theta$, 
we always have $f''_\S(\theta) \theta = 0$ (since $P_{\theta^\bot} \theta = 0$).
Hence, by the symmetry of the second spherical derivative,
$$
\left<f''_\S(\theta)h, \theta\right> = \left<h, f''_\S(\theta) \theta\right> = 0
$$
for any $h \in \R^n$. Hence the image 
$f''_\S(\theta)h$ lies in $\theta^\bot$.

Let us relate the usual second derivative of the classical extension
$u$ as in (7.2) to the spherical derivative $f_\S''$. By (7.3),
$\left<\nabla u(\theta),\theta\right> = 0$, so that
$A = u''(\theta)$ in (7.7). Hence, we get:

\vskip5mm
{\bf Proposition 7.1.} {\sl If $u$ is the $0$-homogeneous extension
of a $C^2$-smooth function $f$ on $\S^{n-1}$ defined in $(7.2)$, then
$$
f_\S''(\theta) = P_{\theta^\bot} u''(\theta) P_{\theta^\bot}, \quad 
\theta \in \S^{n-1}.
$$
}

To argue in the opposite direction from $f$ to $u$, fix a point $x \in \R^n$, 
$x \neq 0$, and put
$$
\theta' = \frac{x+h}{|x+h|}, \quad \theta = \frac{x}{|x|} = r^{-1} x, 
\quad r = |x|,
$$
with $h \in \R^n$ small enough. 
Using Taylor's formula, it is easy to check that
\bee
\theta' - \theta 
 & = &
r^{-1} h - r^{-3}\left<x,h\right> x - r^{-3}\left<x,h\right> h \\
 & & - \ 
\frac{1}{2}\,r^{-3}\,|h|^2\, x + \frac{3}{2}\,
r^{-5}\left<x,h\right>^2 x + O(|h|^3).
\ene
In particular, $|\theta' - \theta| = O(|h|)$. Since 
$\left<\nabla_\S f(\theta),x\right> = r\left<\nabla_\S f(\theta),\theta\right> = 0$,
we get
$$
\left<\nabla_\S f(\theta),\theta' - \theta\right> =
r^{-1} \left<\nabla_\S f(\theta),h\right> - 
r^{-2} \left<\nabla_\S f(\theta),h\right> \left<\theta,h\right> + O(|h|^3).
$$

Recall that the matrix $B = f''_\S(\theta)$ satisfies $Bx = r B\theta = 0$.
Using a shorter expansion
$$
\theta' - \theta = r^{-1} h - r^{-3}\left<x,h\right> x + O(|h|^2),
$$
we also have $B(\theta' - \theta) = r^{-1} Bh + O(|h|^2)$ and therefore
$$
\left<B(\theta' - \theta),\theta' - \theta\right> = 
r^{-2} \left<Bh,h\right> + O(|h|^3).
$$

Hence, by Taylor's expansion (7.5),
\bee
u(x+h) - u(x)
 & = &
f(\theta') - f(\theta) \\
 & & \hskip-15mm = \ 
r^{-1} \left<\nabla_\S f(\theta),h\right> - 
r^{-2} \left<\nabla_\S f(\theta),h\right> \left<\theta,h\right> +
\frac{1}{2}\,r^{-2} \left<Bh,h\right> + o(|h|^2).
\ene
Thus, according to the usual Taylor expansion for the function $u(x)$ 
up to the quadratic term, necessarily 
$\nabla u(x) = r^{-1} \nabla_\S f(\theta)$ and, for all $h \in \R^n$,
\be
\left<u''(x)h,h\right> = - 
2r^{-2} \left<\nabla_\S f(\theta),h\right> \left<\theta,h\right> +
r^{-2} \left<Bh,h\right>.
\en

To give an equivalent matrix description, one may use the symmetrized 
tensor product.

\vskip5mm
{\bf Definition 7.2.} Given vectors $v = (v_1,\dots,v_n), w = (w_1,\dots,w_n)$ 
in $\R^n$, the symmetrized tensor product $v \otimes w$ is
an $n \times n$ symmetric matrix with entries 
$$
(v \otimes w)_{ij} = \frac{1}{2}\,(v_i w_j + w_i v_j), \quad 1 \leq i,j \leq n.
$$

Note that
$$
\left<v,h\right> \left<w,h\right> = \sum_{i,j=1}^n v_i\, w_j h_i h_j = 
\left<(v \otimes w) h,h\right>, \quad h = (h_1,\dots,h_n) \in \R^n.
$$
Using this equality, from (7.8) we obtain:

\vskip5mm
{\bf Proposition 7.3.} {\sl The extension $u(x)$ in $(7.2)$ of a $C^2$-smooth function
$f(\theta)$ on $\S^{n-1}$ has the matrix of second derivatives
$$
u''(x) = r^{-2}\,\big[ f''_\S(\theta) - 2\,\nabla_\S f(\theta) \otimes \theta\big], 
\quad \theta = r^{-1} x, \ x \neq 0.
$$
In particular, on the unit sphere
$$
f''_\S(\theta) =
u''(\theta) + 2\,\nabla_\S f(\theta) \otimes \theta.
$$
}

{\bf Example.} For the linear function $f(\theta) = \left<v,\theta\right>$
with a fixed $v \in \R^n$, we have, by (7.1),
$$
\nabla_\S f(\theta) = P_{\theta^\bot} \nabla f(\theta) = 
v - \left<v,\theta\right> \theta.
$$
Using the extension (7.2) and applying (7.6) or (7.7) together with Proposition 7.3, 
one can show that
$$
f_\S''(\theta) = - \left<v,\theta\right> P_{\theta^\bot} = 
\left<v,\theta\right> [\theta \otimes \theta - I_n].
$$

\vskip10mm
\section{{\bf Spherical Laplacian}}
\setcounter{equation}{0}

\vskip2mm
\noindent
Here we collect basic formulas related to the spherical Laplacian 
with reference to \cite{B-G-L, B-C-G1, B-C-G2} for more details and proofs. 
This operator is defined on the class 
of all $C^2$-smooth functions $f$ on $\S^{n-1}$ by the equality
$$
\Delta_\S f = {\rm Tr}\, f''_\S.
$$

Introduce the ``spherical partial derivatives''
$D_i f(\theta) = \left<\nabla_\S f(\theta),e_i\right>$,
where $e_1,\dots, e_n$ is the canonical basis in $\R^n$, so that
$$
\nabla_\S f(\theta) = \sum_{i=1}^n D_i f(\theta)\, e_i.
$$
As further partial derivatives, one may define the ``second order''
differential operators
$$
D_{ij} f = D_i (D_j f) =
\left<\nabla_\S \left<\nabla_\S f,e_j\right>,e_i\right>,
\quad i,j = 1,\dots,n
$$
(note that the identity $D_{ij} f = D_{ji} f$ is no longer true 
in the entire class $C^2$). Then
$$
\Delta_\S = \sum_{i=1}^n D_{ii}.
$$
In fact, any orthonormal basis in $\R^n$ could be used in
place of the $e_i$'s in the definition of $D_{ii}$, and the
above statement will continue to hold. 

The next statement indicates how one may evaluate the Laplacian of
a given function in terms of usual (Euclidean) derivatives.

\vskip5mm
{\bf Proposition 8.1.} {\sl If $f$ is $C^2$-smooth in an open region $G$ 
containing the unit sphere, then for any $\theta \in \S^{n-1}$,
\be
\Delta_\S f(\theta) = \Delta f(\theta) - (n-1)
\left<\nabla f(\theta),\theta\right> - \left<f''(\theta)\theta,\theta\right>.
\en
}

\vskip2mm
The Laplacian appears naturally in the classical
formula for spherical integration by parts,
\be
\int \left<\nabla_\S f,\nabla_\S g\right> d\sigma_{n-1} = -
\int f \Delta_\S g\,d\sigma_{n-1},
\en
which may be taken as an equivalent definition of $\Delta_\S$. 
It yields the following characterization of the orthogonality to linear 
functions in terms of the spherical gradient.

\vskip5mm
{\bf Proposition 8.2.} {\sl For any $C^1$-smooth function $f$ on $\S^{n-1}$,
$$
\int f(\theta) \theta\,d\sigma_{n-1}(\theta) =
\frac{1}{n-1} \int \nabla_\S f(\theta)\,d\sigma_{n-1}(\theta).
$$
In particular, $f$ is orthogonal to all linear functions in
$L^2(\S^{n-1})$ if and only if the linear forms
$\left<\nabla_\S f(\theta),v\right>$ have $\sigma_{n-1}$-mean
zero for any $v \in \R^n$.
}

\vskip5mm
Applying the rule (7.4), we have another identity
\be
\frac{1}{2}\,\big[\,\Delta_\S(fg) - f \Delta_\S g - g \Delta_\S f\,\big] =
\left<\nabla_\S f,\nabla_\S g\right>
\en
in the class of all $C^2$-smooth functions $f$ and $g$ on the sphere.
For the proof, it is sufficient to multiply this equality by an arbitrary
smooth function $h$ and integrate according to (8.2).

It is well-known (cf.\ e.g.\ \cite{S-W}) that the Hilbert space
$L^2(\S^{n-1},\sigma_{n-1})$ can be decomposed into the sum of
orthogonal linear subspaces $H_d$, $d = 0,1,2,\dots$, consisting of all 
$d$-homogeneous harmonic polynomials (more precisely -- restrictions of 
such polynomials to the sphere). In particular, $H_0$ is the space of constant
functions, $H_1$ is the space of linear functions, $H_2$ consists of quadratic 
harmonics, and so on. Any element $f_d$ of $H_d$ represents an eigenfunction of
the Laplacian with the eigenvalue $- d(n+d-2)$. That is, any function $f$ in 
$L^2(\S^{n-1},\sigma_{n-1})$ admits a unique orthogonal (Fourier) expansion 
in spherical harmonics,
\be
f = \sum_{d=0}^\infty f_d \quad (f_d \in H_d),
\en
and $\Delta_\S f_d = - d(n+d-2)\,f_d$.
Hence, if $f$ is $C^2$-smooth, we get another
representation for the Laplacian,
\be
\Delta_\S f = -\sum_{d=1}^\infty d(n+d-2) f_d.
\en

\vskip5mm
\section{{\bf Semi-group approach to the covariance}}
\setcounter{equation}{0}

\vskip2mm
\noindent
The heat semi-group $P_t = e^{t \Delta_\S}$, $t \geq 0$, 
on the sphere $\S^{n-1}$ has the generator $\Delta_\S$, so that
\be
\frac{d}{dt}\,P_t f = \Delta_\S P_t f, \quad t > 0.
\en
For every $t>0$ and $f \in L^2(\sigma_{n-1})$, the operator $P_t$ assigns
a $C^\infty$-smooth function $P_t f$, which can be defined using the
orthogonal decomposition (8.4) into spherical harmonics by
\be
P_t f = \sum_{d=0}^\infty e^{-d(n+d-2) t}\, f_d.
\en

Let us list a few properties of this semi-group.
To simplify the notations, in this section the gradients $\nabla = \nabla_\S$ 
and the Laplacian $\Delta = \Delta_\S$ are understood in the spherical sense.

\vskip2mm
1) (commutativity) $\Delta P_t f = P_t \,\Delta f$ for any $C^2$-smooth
$f$ on $\S^{n-1}$.

\vskip2mm
2) $\int (P_t f)\,g\,d\sigma_{n-1} = \int f\, P_t g\,d\sigma_{n-1}$.

\vskip2mm
3) $P_t f \rightarrow \int f\,d\sigma_{n-1}$ as $t \rightarrow \infty$.

\vskip2mm
4) (Jensen's inequality) $\Psi(P_t f) \leq P_t \Psi(f)$ for any convex function $\Psi$.

\vskip5mm
There is another important property related to the geometry of 
the sphere known as the curvature-dimension condition (for details we refer 
an interested reader to \cite{B-G-L}). According to (8.3), for the generator 
$L = \Delta$, the carr\'e du champ operator is given by
$$
\Gamma(f,g) = \frac{1}{2}\,\big[L(fg) - f Lg - g Lf\,\big] =
\left<\nabla f,\nabla g\right>.
$$ 
Analogously, one may speak of an iteration of $\Gamma$ 
$$
\Gamma_2(f,g) = \frac{1}{2}\,\big[L\Gamma(f,g) - \Gamma(f,Lg) - \Gamma(g,Lf)\,\big]. 
$$ 

In general, a semi-group is said to satisfy the curvature dimension condition 
$CD(\rho,n)$ with constant $\rho$ if 
$$
\Gamma_2(f) \geq \rho\,\Gamma(f) + \frac{1}{n}\,(Lf)^{2},
$$ 
where $\Gamma(f) = \Gamma(f,f)$. It is a classic result that the heat semi-group 
on $\S^{n-1}$ satisfies the curvature condition $CD(n-2,n-1)$, so that 
$$
\Gamma_2(f) \geq (n-2)\,\Gamma(f)+\frac{1}{n-1}\,(Lf)^{2}.
$$ 
A general characterization commonly known as the Bakry-Emery condition indicates that
$$
CD(\rho,\infty) \iff 
\sqrt{\Gamma\left( P_{t}f\right)}\leq e^{-\rho t}P_{t}\sqrt{\Gamma(f)} \quad 
{\rm for \ all} \ t \geq 0.
$$ 
Hence, the heat semi-group on the sphere satisfies
$$
\sqrt{\Gamma \left(P_t f\right)} \leq e^{-(n-2)t} P_t\sqrt{\Gamma(f)},
$$ 
that is, we have a pointwise bound
\be
|\nabla P_t f|\leq e^{-(n-2)t}P_{t}|\nabla f|.
\en

These properties are sufficient to derive a covariance representation in terms of
the semi-group operators. We start with a variance representation.

\vskip5mm
{\bf Lemma 9.1.} {\sl For any smooth function $f$ on $\S^{n-1}$ $(n \geq 3)$,
\be
\Var_{\sigma_{n-1}}(f) = 
\int_0^\infty \int_{\S^{n-1}} \left<\nabla P_t f,\nabla f\right> d\sigma_{n-1}\,dt.
\en
}

\vskip2mm
{\bf Proof.} We may assume that $f$ is $C^2$-smooth and has mean
$P_\infty f = \int_{\S^{n-1}} f\,d\sigma_{n-1} = 0$. Note that the double
integral in (9.4) is absolutely convergent due to the exponential decay
in (9.3) with respect to $t$. Using $P_0 f = f$ and applying (9.1) together with
properties 3) and 1), we can write
$$
f^2 = -\int_0^\infty \frac{d}{dt}\,(P_t f )^2\,dt = 
-2\int_0^\infty P_t f \, P_t \Delta f\,dt.
$$
Integrating this equality over $\sigma_{n-1}$ and applying 2) together
with Fubini's theorem, we get
\bee
\Var_{\sigma_{n-1}}(f) 
 & = &
-2\int_0^\infty \bigg[\int_{\S^{n-1}} P_{2t}f \, \Delta f\,d\sigma_{n-1}\bigg] dt\\
 & = &
-\int_0^\infty \bigg[\int_{\S^{n-1}} P_t f \,\Delta f\, d\sigma_{n-1}\bigg] dt \\
 & = &
\int_0^\infty \bigg[\int_{\S^{n-1}} 
\left<\nabla P_t f,\nabla f\right>\,d\sigma_{n-1}\bigg] dt,
\ene
where in the last step we employed the formula (8.2) for spherical integration by parts.
\qed

\vskip5mm
{\bf Proof of Theorem 1.2.} In view of the identity
\be
2\,\cov_{\sigma_{n-1}}(f,g) = \Var_{\sigma_{n-1}}(f+g) -
\Var_{\sigma_{n-1}}(f) - \Var_{\sigma_{n-1}}(g),
\en
it makes sense to apply Lemma 8.1 to $f,g$ and $f+g$.
By linearity of the semi-group and spherical gradient, 
\bee
\left<\nabla P_t (f+g),\nabla (f+g)\right>
 & = &
\left<\nabla P_t f,\nabla f\right> +
\left<\nabla P_t g,\nabla g\right> \\
 & & + \ 
\left<\nabla P_t f,\nabla g\right> +
\left<\nabla P_t g,\nabla f\right>.
\ene
Hence, after integration over $d\sigma_{n-1}$ and $dt$, (8.4) yields
\bee
\Var_{\sigma_{n-1}}(f+g) 
 & = &
\Var_{\sigma_{n-1}}(f) + \Var_{\sigma_{n-1}}(g) +
\int_0^\infty \bigg[\int_{\S^{n-1}} 
\left<\nabla P_t f,\nabla g\right> d\sigma_{n-1}\bigg] dt \\
 & & + \ 
\int_0^\infty \bigg[\int_{\S^{n-1}}
\left<\nabla P_t g,\nabla f\right> d\sigma_{n-1}\bigg] dt.
\ene
Here the last inner integral is equal to
\bee
-\int_{\S^{n-1}} f\, \Delta P_t g\, d\sigma_{n-1}
 & = &
-\int_{\S^{n-1}} f\, P_t\Delta g\, d\sigma_{n-1} \, = \,
-\int_{\S^{n-1}} P_t f\, \Delta g\, d\sigma_{n-1}\\ 
 & = &
\int_{\S^{n-1}} \left<\nabla P_t f,\nabla g\right> d\sigma_{n-1}.
\ene
Thus
$$
\Var_{\sigma_{n-1}}(f+g) = \Var_{\sigma_{n-1}}(f) + \Var_{\sigma_{n-1}}(g) +
2\int_0^\infty \bigg[\int_{\S^{n-1}} 
\left<\nabla P_t f,\nabla g\right> d\sigma_{n-1}\bigg] dt.
$$
It remains to apply (8.5) which leads to the desired representation
\be
\cov_{\sigma_{n-1}}(f,g) = 
\int_0^\infty \bigg[\int_{\S^{n-1}} 
\left<\nabla P_t f,\nabla g\right> d\sigma_{n-1}\bigg] dt.
\en
\qed

\vskip5mm
\section{{\bf Comparison of the two representations}}
\setcounter{equation}{0}

\vskip2mm
\noindent
In the Gauss space there is a similar heat semi-group description of the covariance 
functional
\be
\cov_{\gamma_n}(u,v) = \int_0^\infty \bigg[
\int_{\R^n} \left<\nabla\, T_t u,\nabla v\right> d\gamma_n\bigg] dt,
\en
where 
\be
T_t u(x) = \int_{\R^n} u\big(e^{-t}x + \sqrt{1-e^{-2t}} y\big)\,d\gamma_n(y), \quad
x \in \R^n, \ t \geq 0,
\en
denote the Ornstein-Uhlenbeck operators. The identity (10.1) can be obtained
with similar arguments used in the proof of Theorem 1.2. This certainly begs 
the following question: May one start from the Gaussian covariance representation 
(10.1) and obtain the semi-group representation (9.6) of the covariance on the sphere? 
This can be indeed done by applying (10.1) to functions of the form
$u(x) = f(\frac{x}{|x|}) = f(\theta)$, $v(x) = g(\frac{x}{|x|}) = g(\theta)$
by means of orthogonal polynomials. All of the following on Chebyshev-Hermite 
polynomials can be found in the detailed survey by Bogachev 
on the Ornstein-Uhlenbeck operators \cite{Bog}.

In dimension one, the normalized Chebyshev-Hermite polynomial of degree $m \geq 0$ 
is defined as 
$$
H_m(x) = \frac{(-1)^{m}}{\sqrt{m!}}\,e^{x^2/2}\,\frac{d^m}{dx^m}\,e^{-x^2/2}.
$$ 
Similarly, given a multi-index 
$m = (k_1,\dots,k_n)$ with integer components $k_i \geq 0$, the corresponding 
Chebyshev-Hermite polynomial of degree $k = |m| = k_1 + \dots + k_n$ is defined~by
$$
H_m(x) = \prod_{i = 1}^n H_{k_i}(x_i), \quad x = (x_1,\dots,x_n) \in \R^n.
$$
These polynomials constitute an orthonormal basis of $L^2(\gamma_n)$.

The Ornstein-Uhlenbeck operator (10.2) has the following representation 
\be
T_t u = \sum_{k = 0}^\infty e^{-kt}\, I_k(u),
\en
where $I_k$ denotes the projection on the space spanned by the 
$\left( \begin{array}{cc}
k+n-1 \\ 
k 
\end{array}\right)$
multi-index Chebyshev-Hermite polynomials $H_m$ of degree $k$. 

They also serve as eigenfunctions for the generator $L = \Delta - x\cdot\nabla$
of the semi-group $(T_t)_{t \geq 0}$, and more precisely, 
\be
-L H_m(x) = k H_m(x).
\en 
Let us mention the advantage presented 
by an examination of the Ornstein-Uhlenbeck operators from this angle. 
The domain of the generator is given as follows 
$$
\mathcal{D}(L) = \Big\{u \in L^2(\gamma_n): 
\sum_{k=0}^\infty k^2\, ||I_{k}(u)||_{L^{2}(\gamma_{n})}<\infty\Big\}.
$$

Turning to the heat semi-group on $\S^{n-1}$, recall that $f$ in 
$L^2(\sigma_{n-1})$ admits a unique orthogonal decomposition 
$f = \sum_{d = 0}^\infty f_d$ over spherical harmonics of degree $d$. We will 
need the following technical result whose proof we leave
for an interested reader as an exercise.

\vskip5mm
{\bf Lemma 10.1.} {\sl Let $f$ and $g$ be in $L^2(\sigma_{n-1})$.
For the spherical harmonics $f_d$, define 
their extensions $u_d(x)=f_d(\frac{x}{|x|})=f_d(\theta)$ and
let $v(x)=g(\frac{x}{|x|})=g(\theta)$ $(x \in \R^n, x \neq 0)$. Then
\bee
\int_0^\infty \bigg[\int_{\R^n} v \sum_{k = 0}^\infty k e^{-kt} \sum_{d = 0}^\infty
I_k(u_d)\, d\gamma_n \bigg] dt  
 & & \\
 & & \hskip-50mm = \
\sum_{d = 0}^\infty \bigg[d(d+n-2)\,e^{-d(d+n-2)t}
\int_0^\infty \int_{\S^{n-1}} g f_d\,d\sigma_{n-1}\bigg] dt.
\ene
}

Applying this result, let us show that 
$$
\cov_{\sigma_{n-1}}(f,g) = \int_0^\infty \int_{\S^{n-1}}
\left<\nabla_\S P_{t}f,\nabla_\S g \right> d\sigma_{n-1}\,dt.
$$ 
Integrating by parts with respect to the Gaussian measure,
we may rewrite (10.1) as
$$
\cov_{\gamma_n}(u,v) = -\int_0^\infty \int_{\R^n} v L(T_t u)\,d\gamma_n\,dt.
$$
Using a decomposition of $u$ 
into Chebyshev-Hermite polynomials and applying (8.3)-(8.4),
we see that the last double integral is equal to
\bee
-\int_0^\infty \int_{\R^n} v 
L\Big(\sum_{k=0}^\infty e^{-kt}I_{k}(u)\Big)\,d\gamma_n\,dt
 & = &
-\int_0^\infty \int_{\R^n} v\sum_{k=0}^\infty e^{-kt}L(I_{k}(u))\,d\gamma_n\,dt \\
 & = &
\int_0^\infty \int_{\R^n} v\sum_{k=0}^\infty k e^{-kt}\,I_{k}(u)\,d\gamma_n\,dt.
\ene
Next, using the decomposition of $u$ into spherical harmonics and applying
Lemma 10.1 together with the orthogonal decomposition (8.5) for the Laplacian, 
the latter expression is
\bee
\int_0^\infty \int_{\R^n} v \sum_{k = 0}^\infty k e^{-kt}\,
I_k\Big(\sum_{d = 0}^\infty u_d\Big)\,d\gamma_n\, dt
 & = &
\int_0^\infty\int_{\R^n} v \sum_{k = 0}^\infty k e^{-kt} 
\sum_{d=0}^\infty I_k(u_d)\,d\gamma_n\,dt\\
 & = &
\sum_{d = 0}^\infty \bigg[d(d+n-2)\,e^{-d(d+n-2)t}
\int_0^\infty \int_{\S^{n-1}} g f_d\,d\sigma_{n-1}\bigg] dt \\
 & = &
-\int_0^\infty\int_{\S^{n-1}} g \sum_{d=0}^\infty 
e^{-d(d+n-2)t} \Delta_\S f_d\,d\sigma_{n-1}\,dt\\
 & = &
-\int_0^{\infty} \int_{\S^{n-1}} g\, \Delta_\S \sum_{d=0}^\infty 
e^{-d(d+n-2)t} f_d\,d\sigma_{n-1}\,dt\\
 & = &
-\int_0^\infty\int_{\S^{n-1}} g \,\Delta_\S P_t f\,d\sigma_{n-1}\,dt\\
 & = &
\int_0^\infty\int_{\S^{n-1}} 
\left<\nabla_\S P_t f,\nabla_\S g\right> d\sigma_{n-1}\,dt.
\ene

\vskip5mm
\section{{\bf Applications to spherical concentration}}
\setcounter{equation}{0}

\vskip2mm
\noindent
We now return to Theorem 1.1 and develop several applications of the 
covariance identity
\be
{\rm cov}_{\sigma_{n-1}}(f,g) = c_n \int_{\S^{n-1}} \int_{\S^{n-1}} 
\left<\nabla_\S f(x),\nabla_\S g(y)\right> d\mu_n(x,y).
\en
All of them are entirely based on the property that the probability measure 
$\mu_n$ on $\S^{n-1} \times \S^{n-1}$ (explicitly described in Theorem 4.1)
has marginals $\sigma_{n-1}$, while the constants satisfy
$$
\frac{1}{n-1} < c_n < \frac{1}{n-2} \quad (n \geq 3).
$$

As a particular case, (11.1) implies the Poincar\'e-type inequality
$$
{\rm Var}_{\sigma_{n-1}}(f) \leq c_n \int_{\S^{n-1}}
|\nabla_\S f|^2\,d\sigma_{n-1}
$$
with an asymptotically correct constant 
(the optimal one is $\frac{1}{n-1}$). More generally, a similar bound can be
given in terms of $L^p$-norms of gradients
$$
\|\nabla_\S f\|_p = \Big(\int_{\S^{n-1}}
|\nabla_\S f|^p\,d\sigma_{n-1}\Big)^{1/p}, \quad 1 \leq p \leq \infty,
$$
where the case $p=\infty$ corresponds to the maximum (for smooth $f$).

\vskip5mm
{\bf Corollary 11.1.} {\sl Let $p,q \geq 1$, $\frac{1}{p} + \frac{1}{q} = 1$,
For all smooth functions $f$ and $g$ on $\S^{n-1}$,
\be
|{\rm cov}_{\sigma_{n-1}}(f,g)| \leq c_n\,\|\nabla_\S f\|_p\, \|\nabla_\S g\|_q.
\en
In particular, if $f$ has a Lipschitz semi-norm $\|f\|_{\rm Lip} \leq 1$, then
\be
|{\rm cov}_{\sigma_{n-1}}(f,g)| \leq c_n \int_{\S^{n-1}}
|\nabla_\S g|\, d\sigma_{n-1}.
\en
}

\vskip2mm
{\bf Proof.} By H\"older's inequality, the absolute value of the double integral 
in (10.1) does not exceed
\bee
\int_{\S^{n-1}} \int_{\S^{n-1}} |\nabla_\S f(x)|\,|\nabla_\S g(y)|\, d\mu_n(x,y) 
 & \leq & \\ 
 & & \hskip-60mm
\Big(\int_{\S^{n-1}} \int_{\S^{n-1}} |\nabla_\S f(x)|^p\, d\mu_n(x,y)\Big)^{1/p}\,
\Big(\int_{\S^{n-1}} \int_{\S^{n-1}} |\nabla_\S g(y)|^q\, d\mu_n(x,y)\Big)^{1/q} \, = \, \\
 & & \hskip-60mm
\Big(\int_{\S^{n-1}} |\nabla_\S f(x)|^p\, d\sigma_{n-1}(x)\Big)^{1/p}\,
\Big(\int_{\S^{n-1}} |\nabla_\S g(y)|^q\, d\sigma_{n-1}(y)\Big)^{1/q}\, = \,
\|\nabla_\S f\|_p\, \|\nabla_\S g\|_q.
\ene
The inequality (11.3) corresponds to (11.2) with $p = \infty$, $q = 1$.
\qed

\vskip5mm
The representation (11.1) may also be used to recover the spherical concentration.

\vskip5mm
{\bf Corollary 11.2.} {\sl For any smooth function $f$ on $\S^{n-1}$ with
$\sigma_{n-1}$-mean zero,
\be
\int_{\S^{n-1}} e^f\,d\sigma_{n-1} \leq \int_{\S^{n-1}}\int_{\S^{n-1}} 
\exp\big\{c_n \left<\nabla_\S f(x),\nabla_\S g(y)\right>\big\}\, d\mu_n(x,y).
\en
In particular,
\be
\int_{\S^{n-1}} e^f\,d\sigma_{n-1} \leq \int_{\S^{n-1}}
e^{c_n |\nabla_\S f|^2}\,d\sigma_{n-1}.
\en
}

\vskip2mm
{\bf Proof.} The argument involves the entropy functional
$$
\Ent(\xi) = \E\,\xi \log \xi  - \E\,\xi \log \E\,\xi.
$$
It is well-defined and non-negative for every random variable
$\xi \geq 0$, and is known to admit the variational description
$$
\Ent(\xi) = \sup \E\, \xi \eta,
$$
where the supremum is taken over all random variables $\eta$ such that
$\E e^\eta \leq 1$. In particular,
\be
\E e^\eta = 1 \ \Longrightarrow \ \E\,\xi\eta \leq \Ent(\xi).
\en

We consider the expectations (integrals) and entropy on the probability spaces
$(\S^{n-1}, \sigma_{n-1})$ and $(\S^{n-1} \times \S^{n-1}, \mu_n)$ and 
write correspondingly $\E_{\sigma_{n-1}} \xi$, $\Ent_{\sigma_{n-1}}(\xi)$, and 
$\E_{\mu_n} \xi$, $\Ent_{\mu_n}(\xi)$. 

Define the constant $\beta$ as the logarithm 
of the right-hand side in (11.4), so that, by (11.6),
\be
\E_{\mu_n} \Big[\big(c_n \left<\nabla_\S f(x),\nabla_\S g(y)\right> - \beta\big) \xi\Big]
 \leq \Ent_{\mu_n}(\xi)
\en
for any bounded measurable function $\xi(x,y)$ on $\S^{n-1} \times \S^{n-1}$. 
Choosing $\xi(x,y) = g(y) = e^{f(y)}$, we have 
$\Ent_{\mu_n}(\xi) = \Ent_{\sigma_{n-1}}(g)$, and (11.1) together with (11.7) give
\bee
\E_{\sigma_{n-1}} f e^f - \beta\, \E_{\sigma_{n-1}} e^f
 & = &
{\rm cov}_{\sigma_{n-1}}(f,e^f) - \beta\, \E_{\sigma_{n-1}} e^f \\
 & = & 
c_n\, \E_{\mu_n} \left<\nabla_\S f(x),\nabla_\S g(y)\right> e^{f(y)} - 
\beta\, \E_{\mu_n} e^f \\ 
 & \leq &
\Ent_{\mu_n}(e^f) = \Ent_{\sigma_{n-1}}(e^f) \\
 & = &
\E_{\sigma_{n-1}}(f e^f) - \E_{\sigma_{n-1}} e^f \log \E_{\sigma_{n-1}} e^f.
\ene
Hence
$\beta \geq \log \E_{\sigma_{n-1}} e^f$, which is the relation (11.4).
As for (11.5), it follows from (11.4) using
$$
\left<\nabla_\S f(x),\nabla_\S g(y)\right> \leq 
\frac{1}{2}\,|\nabla_\S f(x)|^2 + \frac{1}{2}\,|\nabla_\S f(y)|^2
$$ 
and applying the Cauchy inequality.
\qed

\vskip5mm
{\bf Proof of Theorem 1.3.} The argument is based on the relation (11.2),
where the smoothness condition on $f$ may be removed, just keeping
the assumption $\|f\|_{\rm Lip} \leq 1$. We may further assume that
$f$ has $\sigma_{n-1}$-mean zero. 
Consider $g = T(f)$, where $T$ is a non-decreasing (piecewise) differentiable 
function. Then from (11.2) we get
$$
\E_{\sigma_{n-1}}\, f T(f) \leq c_n\, \E_{\sigma_{n-1}}\, T'(f).
$$
Without loss of generality, assume that $f$ has a continuous density
$p$ under the measure $\sigma_{n-1}$. Let $F(x) = \sigma_{n-1}\{f \leq x\}$ 
denote the distribution function of $f$. Choosing $T(x) = \min((x-r)^+,\ep)$
with parameters $r>0$ and $\ep>0$, the above inequality becomes
$$
\int_r^{r+\ep} x(x-r)\,dF(x) + \ep \int_{r + \ep}^\infty x\,dF(x) \leq
c_n\,(F(r+\ep) - F(r)).
$$
Dividing by $\ep$ and letting $\ep$ tend to zero, we obtain that for all 
$r > 0$,
$$
V(r) \equiv \int_r^\infty xp(x)\,dx = \int_r^\infty x\,dF(x) \leq c_n p(r).
$$
Hence, the function $V$ satisfies the differential inequality
$V(r) \leq - c_n V'(r)/r$, that is,
$$
(\log V(r))' \leq (-r^2/2c_n)'.
$$
This is the same as saying that $\log V(r) + r^2/2c_n$
is non-increasing, and so is the function
$V(r)\,\exp(r^2/2c_n)$. Since $V(0) = \E\,f^+$, we get
$$
\E_{\sigma_{n-1}}\,f^+ \geq \exp(c_n r^2/2)\int_r^\infty xp(x)\,dx \geq
\exp(r^2/2c_n)\,r(1-F(r)).
$$
Thus,
$$
\sigma_{n-1}\{f \geq r\} \leq \frac{1}{r}\,e^{-r^2/2c_n}\, 
\E_{\sigma_{n-1}}\,f^+.
$$
A similar inequality holds when replacing $f$ with $-f$, that is,
$$
\sigma_{n-1}\{f \leq -r\} \leq \frac{1}{r}\,e^{-r^2/2c_n}\, 
\E_{\sigma_{n-1}}\,f^-.
$$
Summing the two bounds, we arrive at (1.8).
\qed


\vskip10mm
\section{{\bf Second order covariance identity in Gauss space}}
\setcounter{equation}{0}

\vskip2mm
\noindent
For short, we write $\E_{\gamma_n}$ to denote the expectation, that is, the
integral with respect to the Gaussian measure $\gamma_n$ on $\R^n$.

Applying the Gaussian covariance identity (1.1)
to partial derivatives of $u$ and $v$ of the first order, 
we are led to the second covariance identity
\be
{\rm cov}_{\gamma_n}(u,v) = 
\left<\E_{\gamma_n} \nabla u,\E_{\gamma_n} \nabla v\right> + 
\int_{\R^n}\!\int_{\R^n} \left<u''(x),v''(y)\right> d\kappa_n(x,y),
\en
where the integrand represents the inner product of square symmetric matrices, 
which is the trace ${\rm Tr}(u''(x)\, v''(y))$ 
of the product of the two Hessians. Here, $2\kappa_n$ is a probability measure 
on $\R^n \times \R^n$ which may be defined by means of the mixture
\be
\kappa_n = \int_0^1 (1-t)\, {\sf L}\big(X,tX + sZ\big)\,dt
\en
with convention that $s = \sqrt{1-t^2}$. This shows that $2\kappa_n$ has $\gamma_n$ 
as marginals, similarly to the mixing measure $\pi_n$ in the first order 
covariance representation for the Gaussian measure $\gamma_n$. 

In an equivalent form, the identity (12.1) was emphasized in \cite{H-P-S}.
It may be stated in the class of all $C^2$-smooth
functions $u$ and $v$ on $\R^n \setminus \{0\}$ such that
$$
\E_{\gamma_n} \|u''\|_{\rm HS}^2 < \infty, \quad
\E_{\gamma_n} \|v''\|_{\rm HS}^2 < \infty.
$$
Then the expectations 
$$
\E_{\gamma_n} (u^2 + v^2) \quad {\rm and } \quad
\E_{\gamma_n} (|\nabla u|^2 + |\nabla v|^2)
$$ 
are also finite, by the Poincar\'e-type inequality (2.3) for $\gamma_n$, 
so that (12.1) makes sense.

In addition, (12.2) implies that $\kappa_n$ is absolutely continuous 
with respect to the Lebesgue measure on $\R^n \times \R^n$ and has 
some density
$$
q_n(x,y) = \frac{d\kappa_n(x,y)}{dx\,dy}.
$$
To write it explicitly, recall that for any bounded Borel measurable function 
$h:\R^n \times \R^n \rightarrow \R$, changing the variable $tx + sz = y$,  
or $z = \frac{y - tx}{s}$, we have
\bee
\E\, h(X,tX + sZ) 
 & = & 
\int_{\R^n}\!\int_{\R^n}
h(x,tx + sz)\,\varphi_n(x)\varphi_n(z)\ dx\,dz \\
 & = & 
s^{-n} \int_{\R^n}\!\int_{\R^n} 
h(x,y)\,\varphi_n(x)\varphi_n\Big(\frac{y-tx}{s}\Big)\ dx\,dy,
\ene
and this means that
$$
q_n(x,y) = \int_0^1 (1-t) s^{-n}\,
\varphi_n(x)\varphi_n\bigg(\frac{y-tx}{s}\bigg)\,dt.
$$
But
$$
\varphi_n(x)\varphi_n\bigg(\frac{y-tx}{s}\bigg) =
\frac{1}{(2\pi)^n} \ \exp\bigg[-
\frac{|x|^2 + |y|^2 -2t\left<x,y\right>}{2s^2}\bigg].
$$
Hence, one may complement (12.1) with the following assertion,
where we write $s = \sqrt{1-t^2}$.

\vskip5mm
{\bf Lemma 12.1.} {\sl The mixing measure $\kappa_n$ 
in the Gaussian covariance representation $(12.1)$ has density
$$
q_n(x,y) = \frac{1}{(2\pi)^n} \,
\int_0^1 (1-t) s^{-n} \,\exp\bigg[-
\frac{|x|^2 + |y|^2 -2t\left<x,y\right>}{2s^2}\bigg]\,dt, \quad
x,y \in \R^n.
$$
}

\vskip5mm
Putting $x = r\theta$, 
$y = r'\theta'$ ($r,r'>0$, $\theta, \theta' \in \S^{n-1}$), this density
may be written in polar coordinates as
\be
q_n(r\theta,r'\theta') = \frac{1}{(2\pi)^n} \,
\int_0^1 (1-t) s^{-n} \,\exp\bigg[-
\frac{r^2 + r'^2 -2rr' t\left<\theta,\theta'\right>}{2s^2}\bigg]\,dt.
\en


\vskip10mm
\section{{\bf Second order covariance identities on the sphere}}
\setcounter{equation}{0}

\vskip2mm
\noindent
We now develop a spherical variant of the second order covariance identity
(12.1). With any $C^2$-smooth $f$ on $\S^{n-1}$ we associate the function 
$$
u(x) = f(r^{-1} x) = f(\theta), \quad {\rm where} \ \ r = |x|, \ \theta = r^{-1} x,
$$
which is defined and $C^2$-smooth in $\R^n \setminus\{0\}$. By (7.3)
and Proposition 7.3, it has derivatives
\bee
\nabla u(x) 
 & = & r^{-1}\, \nabla_\S f(\theta), \\
u''(x) 
 & = &
r^{-2}\,Df(\theta) \, = \, r^{-2}\,
\big(f''_\S(\theta) - 2\,\nabla_\S f(\theta) \otimes \theta\big)
\ene
(cf. Definition 7.2 for the symmetrized tensor product). Hence, as in Section 4,
$$
\int \nabla u\,d\gamma_n = \E\,\frac{1}{|Z|}\,
\int_{\S^{n-1}} \nabla_\S f\,d\sigma_{n-1}
$$
and
$$
\int \|u''\|_{\rm HS}^2\,d\gamma_n = \E\,\frac{1}{|Z|^4}\,
\int_{\S^{n-1}} 
\|f''_\S(\theta) - 2\,\nabla_\S f(\theta) \otimes \theta\|_{\rm HS}^2\,d\sigma_{n-1},
$$
where $Z$ is a standard normal random vector in $\R^n$.
By the smoothness assumption, the last integrand is bounded (so, the
last integral is finite), while the expectation is finite for $n \geq 5$, only,
in which case, by Lemma 3.1,
$$
\E\,\frac{1}{|Z|^4} = \frac{1}{n-4}\,\E\,\frac{1}{|Z|^2} = \frac{1}{(n-2)(n-4)}.
$$

Now, let $g$ be another $C^2$-smooth function on $\S^{n-1}$.
Define $v(y) = g(\theta')$ for $y \in \R^n \setminus \{0\}$,
$\theta' = \frac{y}{r'}$, $r' = |y|$, so that
$$
\nabla v(y) = \frac{1}{r'}\, \nabla_\S g(\theta') \quad {\rm and} \quad
v''(y) = \frac{1}{r'^2}\,D g(\theta').
$$
To simplify the resulting formulas, let us assume that $f$ and $g$ are orthogonal 
to linear functions in $L^2(\sigma_{n-1})$. This is equivalent to the similar 
property of $u$ and $v$ in $L^2(\gamma_n)$, which is the same as saying that 
$\nabla u$ and $\nabla v$ have $\gamma_n$-mean zero (compare with Proposition 8.2). 
We may now apply (12.1) together with Lemma 12.1, which give
\be
{\rm cov}_{\sigma_{n-1}}(f,g) = {\rm cov}_{\gamma_n}(u,v)
 \, = \,
\int_{\R^n}\!\int_{\R^n} \left<u''(x),v''(y)\right> q_n(x,y)\,dx\,dy.
\en

To continue, let us integrate in (13.1) in polar coordinates separately along
the $x$ and $y$ variables in the same manner as we did in Section 4.
Then, according to (12.3), we obtain 
\be
{\rm cov}_{\sigma_{n-1}}(f,g) = \int_{\S^{n-1}}\int_{\S^{n-1}} 
\left<f_\S''(\theta),g_S''(\theta')\right>
\psi_n(\theta,\theta')\,d\sigma_{n-1}(\theta)\, d\sigma_{n-1}(\theta')
\en
with $\psi_n(\theta,\theta')$ representable as
$$
\frac{1}{2^{n - 2}\, \Gamma(\frac{n}{2})^2}
\int_0^1 (1-t) s^{-n} \bigg[\int_0^\infty\!\int_0^\infty 
\exp\bigg[-\frac{r^2 + r'^2 -2rr' t\left<\theta,\theta'\right>}{2s^2}\bigg]
(r r')^{n-3} \,dr\, dr'\bigg] dt.
$$
Changing the variable, this expression is simplified to
\be
\frac{1}{2^{n-2}\, \Gamma(\frac{n}{2})^2}
\int_0^1 (1-t) s^{n-4} \, \bigg[\int_0^\infty\!\int_0^\infty 
\exp\bigg[-\frac{r^2 + r'^2 -2rr' t\left<\theta,\theta'\right>}{2}\bigg]
(r r')^{n-3} \,dr\, dr'\bigg] dt.
\en

Let us see that $\psi_n$ is integrable over the product measure
$\sigma_{n-1} \otimes \sigma_{n-1}$ and serves as density of some 
finite positive measure on $\S^{n-1} \times \S^{n-1}$ with total mass
$$
c_n = \kappa_n(\S^{n-1} \times \S^{n-1}) = 
\int_{\S^{n-1}}\int_{\S^{n-1}} \psi_n(\theta,\theta')\,
d\sigma_{n-1}(\theta)\,d\sigma_{n-1}(\theta').
$$
Note that, for the finiteness of (13.3), it is necessary that the integral
$\int_0^1 (1-t) s^{n-4}\,dt$ be convergent, that is, $n \geq 3$.
Repeating integration in polar coordinates as before, we have
$$
\E \int_0^1 \frac{1}{|X|^2\,|tX + sZ|^2}\,dt
=
\int_{\R^n}\!\int_{\R^n} \frac{1}{|x|^2\,|y|^2}\ q_n(x,y)\,dx\, dy = c_n,
$$
where $X$ and $Z$ are independent standard normal random vectors in 
$\R^n$. But, for any fixed $0 < t < 1$, by Cauchy's inequality,
\bee
\E\, \frac{1}{|X|^2\,|tX + sZ|^2} 
 & \leq &
\big(\E\,|X|^{-4}\big)^{1/2}\, \big(\E\,|tX + sZ|^{-4}\big)^{1/2} \\
 & = &
\E\,|X|^{-4} \, = \, \frac{1}{(n-2)(n-4)}.
\ene
Thus, $c_n < \frac{1}{(n-2)(n-4)}$ and $\kappa_n = c_n \mu_n$ for some 
probability measure $\mu_n$ on $\S^{n-1} \times \S^{n-1}$.

We have also a natural lower bound on this constant. Indeed, 
for a fixed value of $X$,
$$
\E\,|tX + sZ|^2 = \E\,(t^2 |X|^2 + s^2|Z|^2 + 2ts \left<X,Z\right>) =
t^2 |X|^2 + s^2 n.
$$
Hence, using $\E\,|X|^4 = n^2 + 2n$ and $t^2 + s^2 = 1$, we get
\bee
\E\,|X|^2 \,|tX + sZ|^2 
 & = &
\E_X\,|X|^2\,\Big(\E_Z\, |tX + sZ|^2\Big) \\
 & = &
\E_X\,|X|^2\,(t^2 |X|^2 + s^2 n) \, = \, n^2 + 2t^2n \, \leq \, n^2 + 2n.
\ene
Hence, by Jensen's inequality, using the convexity of the function
$1/r$ in $r>0$,
$$
c_n \, = \, \int_0^1 \E\, \frac{1}{|X|^2\,|tX + sZ|^2}\,dt \, \geq \,
\int_0^1 \frac{1}{\E\,|X|^2 |tX + sZ|^2}\,dt \,\geq \, \frac{1}{n(n+2)}.
$$

Note also that 
$\psi_n(\theta,\theta') = \psi_n(\left<\theta,\theta'\right>)$ 
depends on $\theta$ and $\theta'$ via the inner product only.
This implies that the integral
$$
\int_{S^{n-1}} \psi_n(\theta,\theta')\,d\sigma_{n-1}(\theta') =
\int_{S^{n-1}} \psi_n(\left<\theta,\theta'\right>)\,
d\sigma_{n-1}(\theta')
$$
does not depend on $\theta$ and is therefore equal to $c_n$. 
In other words, the marginal distributions
of $\mu_n$ represent spherically invariant measures,
hence coincide with $\sigma_{n-1}$. 

As a result, we obtain the following spherical analog of (12.1),
using the second order linear matrix-valued differential operator
\be
Df(x) = f''_\S(x) - 2\,\nabla_\S f(x) \otimes x, \quad x \in \S^{n-1}.
\en

\vskip5mm
{\bf Theorem 13.1.} {\sl On $\S^{n-1} \times \S^{n-1}$ $(n \geq 5)$ there 
exist a probability measure $\mu_n$ with marginals $\sigma_{n-1}$
and a constant $c_n>0$ such that
\be
{\rm cov}_{\sigma_{n-1}}(f,g) = c_n \int_{\S^{n-1}} \int_{\S^{n-1}} 
\left<Df(x),Dg(y)\right> d\mu_n(x,y)
\en
for all $C^2$-smooth $f,g$ on $\S^{n-1}$ orthogonal to all linear 
functions in $L^2(\sigma_{n-1})$. Moreover, $\mu_n$ has density with respect to 
the product measure $\sigma_{n-1} \otimes \sigma_{n-1}$ of the form
$\psi_n(\left<x,y\right>)$ for the positive function
$\psi_n$ given in $(13.3)$, and
\be
\frac{1}{n(n+2)} < c_n < \frac{1}{(n-2)(n-4)}.
\en
}

\vskip2mm
One may also develop a heat semi-group variant of the identity (13.5)
in a full analogy between the covariance representations (1.4) and (1.7).
As a consequence of (1.7) we have the following identity where the orthogonality
hypothesis is not needed.

\vskip5mm
{\bf Theorem 13.2.} {\sl For all $C^2$-smooth $f,g$ on $\S^{n-1}$ $(n \geq 3)$,
\be
\cov_{\sigma_{n-1}}(f,g) = 
\int_0^\infty\bigg[\int_{\S^{n-1}} t\,
\Delta_\S P_t f \, \Delta_\S g\, d\sigma_{n-1}\bigg] dt.
\en
}

\vskip2mm
{\bf Proof.} Applying the pointwise bound (9.3) together with Cauchy's inequality,
we have that, for all $t \geq 0$,
\bee
t\, |\left<\nabla_\S P_t f, \nabla_\S g\right>| 
 & \leq &
t\,|\nabla P_t f| \,|\nabla_\S g| \\
 & \leq & t e^{-(n-2)t} P_t |\nabla f| \,|\nabla_\S g|.
\ene
Here, the latter expression is vanishing at zero, and the same is true at $t = \infty$,
since $P_t |\nabla f| \rightarrow \int_{\S^{n-1}}|\nabla f|\,d\sigma_{n-1}$
as $t \rightarrow \infty$. This justifies integration by parts in the
representation (1.4):
\bee
\cov_{\sigma_{n-1}}(f,g) 
 & = &
\int_{\S^{n-1}} \bigg[\int_0^\infty
\left<\nabla_\S P_t f,\nabla_\S g\right> dt \bigg] d\sigma_{n-1} \nonumber \\
 & = &
- \ \int_{\S^{n-1}} \bigg[\int_0^\infty t\, \frac{d}{dt}
\left<\nabla_\S P_t f,\nabla_\S g\right> dt\bigg] d\sigma_{n-1} \nonumber \\
 & = &
- \ \int_0^\infty t\, \bigg[\int_{\S^{n-1}} \frac{d}{dt}
\left<\nabla_\S P_t f,\nabla_\S g\right> d\sigma_{n-1}\bigg] dt.
\ene
Moreover, differentiating under the scalar product sign according to the
semi-group identity (9.1) and applying the formula (8.2) on the spherical 
integration by parts, we see that the last inner integral is equal to
$$
\int_{\S^{n-1}}
\left<\nabla_\S \Delta_\S P_t f,\nabla_\S g\right> d\sigma_{n-1} = -
\int_{\S^{n-1}} \Delta_\S P_t f \, \Delta_\S g\, d\sigma_{n-1}.
$$
\qed

\vskip2mm
{\bf Remark 13.3.} The last step in this derivation can be changed
to yield an equivalent representation in place of (12.7), namely
\be
\cov_{\sigma_{n-1}}(f,g) = \int_0^\infty\bigg[\int_{\S^{n-1}} t
P_t f \, \Delta_\S^2 g\, d\sigma_{n-1}\bigg] dt,
\en
where $\Delta_\S^2 f = \Delta_\S\, \Delta_\S f$ is the square of the spherical
Laplacian.


\vskip10mm
\section{{\bf Upper bounds on covariance of order $1/n^2$}}
\setcounter{equation}{0}

\vskip2mm
\noindent
Note that the constants (13.6) are of order $1/n^2$ in a big contrast 
with constants in the covariance representation of Theorem 1.1 which have
the rate $1/n$ (for growing $n$). It is therefore reasonable to derive
from (13.5) a corresponding analog of Corollary 11.1 in terms of $L^p$-norms
of spherical derivatives and Hessians 
\be
\|\nabla_\S f\|_p = 
\Big(\int_{\S^{n-1}} |\nabla_\S f|^p\,d\sigma_{n-1}\Big)^{1/p}, \quad
\|f''_\S\|_p = 
\Big(\int_{\S^{n-1}} \|f_\S''\|_{\rm HS}^p\,d\sigma_{n-1}\Big)^{1/p}
\en
in the class of functions as in Theorem 12.1.

Since the inner product of $n \times n$ matrices satisfies
$|\left<A,B\right>| \leq \|A\|_{\rm HS}\,\|B\|_{\rm HS}$
(Cauchy's inequality), we have, for all $x,y \in \S^{n-1}$, 
\be
|\left<Df(x),Dg(y)\right>| \leq \|D f(x)\|_{\rm HS}\,\|D f(y)\|_{\rm HS}.
\en

Now, by the triangle inequality for the Hilbert-Schmidt norm, from (13.4) 
it follows that
$$
\|D f(x)\|_{\rm HS} \, \leq \, \|f_S''(x)\|_{\rm HS} + 
2\, \|\nabla_\S f(x) \otimes x\|_{\rm HS}.
$$
In turn, according to Definition 7.2, the symmetrized tensor product
$\nabla_\S f(x) \otimes x$ has entries
$$
(\nabla_\S f(x) \otimes x)_{ij} = \frac{1}{2}\,\big(x_i D_j f(x) + x_j D_i f(x)\big),
\quad 1 \leq i,j \leq n,
$$
so that
$$
(\nabla_\S f(x) \otimes x)_{ij}^2 \leq \frac{1}{2}\,
\Big(x_i^2 (D_j f(x))^2 + x_j^2 (D_i f(x))^2\Big)
$$
and
$$
\|\nabla_\S f(x) \otimes x)\|_{\rm HS}^2 \, \leq \, \sum_{i=1}^n (D_i f(x))^2 \, = \,
|\nabla_\S f(x)|^2
$$
for all $x \in \S^{n-1}$. Thus,
\be
\|D f(x)\|_{\rm HS} \, \leq \, \|f_\S''(x)\|_{\rm HS} + 
2\, |\nabla_\S f(x)|,
\en
and (14.2) gives
$$
|\left<Df(x),Dg(y)\right>| \leq 
\Big(\|f_S''(x)\|_{\rm HS} + 2\, |\nabla_\S f(x)|\Big)\,
\Big(\|g_S''(y)\|_{\rm HS} + 2\, |\nabla_\S g(y)|\Big).
$$
As a consequence from (13.5), ${\rm cov}_{\sigma_{n-1}}(f,g)$ is bounded in absolute 
value by
\be
c_n \int_{\S^{n-1}} \int_{\S^{n-1}} 
\Big(\|f_S''(x)\|_{\rm HS} + 2\, |\nabla_\S f(x)|\Big)\,
\Big(\|g_S''(y)\|_{\rm HS} + 2\, |\nabla_\S g(y)|\Big) d\mu_n(x,y).
\en

Moreover, since the marginals of the probability measure $\mu_n$ in Theorem 13.1
coincides with $\sigma_{n-1}$, the $L^p$-norm of the expression in the first 
bracket in (14.4) with respect to $\mu_n$ and $\sigma_{n-1}$ are equal to each other
and does not exceed $\|f''_\S\|_p + 2\,\|\nabla_\S f\|_p$.
A similar conclusion is true about the $L^q$-norm of the expression in the second
bracket, which is needed when $q$ is the conjugate power.
Applying H\"older's inequality, we may conclude.

\vskip5mm
{\bf Corollary 14.1.} {\sl Let $p,q \geq 1$, $\frac{1}{p} + \frac{1}{q} = 1$.
Given $C^2$-smooth functions $f,g$ on $\S^{n-1}$ $(n \geq 5)$ which are orthogonal 
to all linear functions in $L^2(\sigma_{n-1})$, we have
\be
|{\rm cov}_{\sigma_{n-1}}(f,g)| \leq \frac{1}{(n-2)(n-4)}\,
\Big(\|f''_\S\|_p + 2\,\|\nabla_\S f\|_p\Big) \,
\Big(\|g''_\S\|_q + 2\,\|\nabla_\S g\|_q\Big),
\en
where the $L^p$ and $L^q$ norms are defined in $(14.1)$.
}

\vskip5mm
If these norms are of order at most 1, the covariance $f$ and $g$ will be
therefore of order at most $1/n^2$ for the growing dimension $n$. As easy to see, 
the inequality such as (14.5) remains to hold for $2 \leq n \leq 4$ with
an absolute constant in front of the norms.

In the case $p=q=2$, (14.5) can be simplified by removing the norms for the 
spherical gradients. Indeed, as was shown in \cite{B-C-G1},
$\|\nabla_\S f\|_2 \leq \|f''_\S\|_2$
for any $C^2$-smooth function $f$ on $\S^{n-1}$ 
(cf. also \cite{B-C-G2}, Proposition 10.1.2). In fact, if $f$ is orthogonal 
to all linear functions in $L^2(\sigma_{n-1})$, a much stronger inequality
holds true, namely
$$
\|\nabla_\S f\|_2^2 \leq \frac{1}{n+2}\,\|f''_\S\|_2^2.
$$
Applying this in (14.5) with $f=g$, we obtain a second order Poincar\'e-type
inequality
$$
\Var_{\sigma_{n-1}}(f) \leq \frac{c}{n^2}\,\|f''_\S\|_2^2
$$
up to an absolute constant $c$. 

In order to get a more precise relation in this particular case, one may actually 
employ the following identity for the $L^2$-norm of the Hessian: 
Given a $C^4$-smooth function $f$ on $\S^{n-1}$,
$$
\|f''_\S\|_2^2 = 
\int_{\S^{n-1}} f\,(\Delta_\S^2 f + (n-2)\, \Delta_\S f)\,d\sigma_{n-1},
$$
where we recall that $\Delta_\S^2 f = \Delta_\S\, \Delta_\S f$. 
Using the representation (8.5) of this operator in terms of spherical
harmonics, it follows from this identity that
\be
\Var_{\sigma_{n-1}}(f) \leq \frac{1}{2n (n+2)}\,\|f''_\S\|_2^2
\en
for any $C^2$-smooth function $f$ on $\S^{n-1}$ which is orthogonal 
to all linear functions in $L^2(\sigma_{n-1})$. Here, an equality is attained
for quadratic harmonics, cf. \cite{B-C-G1}.

Another approach to (14.6)
can be based on the semi-group covariance representation (13.8), cf. Remark 13.3.


\vskip10mm
\section{{\bf Second order concentration on the sphere}}
\setcounter{equation}{0}

\vskip2mm
\noindent
Here we employ Theorem 13.1 to strengthen Theorem 1.3 with respect to the growing 
dimension $n$ for a certain class of functions.
The next assertion is analogous to Corollary 10.2; however, we will use 
the constants $c_n$ from (13.6) which are of the order $1/n^2$.

\vskip5mm
{\bf Corollary 15.1.} {\sl Given a $C^2$-smooth function $f$ on the sphere
$\S^{n-1}$ $(n \geq 5)$ which is orthogonal to all affine functions 
in $L^2(\sigma_{n-1})$, we have
\be
\int_{\S^{n-1}} e^f\,d\sigma_{n-1} \leq \int_{\S^{n-1}} 
\exp\Big\{\frac{1}{(n-2)(n-4)}\,
(2\,\|f''_\S\|_{\rm HS}^2 + 8\,|\nabla_\S f|^2)\Big\}\,d\sigma_{n-1}.
\en
}

\vskip2mm
A similar property was proved in \cite{B-C-G1} under the following assumptions:

\vskip2mm
$a)$ \ $f$ is orthogonal to all affine functions in $L^2(\sigma_{n-1})$;

\vskip1mm
$b)$ \ $\|f''_\S\| \leq 1$ pointwise on the sphere where
$\|f''_\S\|$ denotes the operator norm;

\vskip1mm
$c)$ \ $\int \|f''_\S\|_{\rm HS}^2\,d\sigma_{n-1} \leq b$.

\vskip2mm
\noindent
Then it was shown that
\be
\int_{\S^{n-1}} \exp\Big\{\frac{n-1}{2 (1+b)}\,|f|\Big\}\,d\sigma_{n-1} \leq 2.
\en
This implies that $f$ is of the order at most $1/n$, and we have a deviation
inequality
$$
\sigma_{n-1}\big\{(n-1) |f| \geq r\big\} \leq 2 e^{-r/2(1+b)}, \quad r \geq 2.
$$

A similar conclusion can also be made on the basis of (15.1) if the conditions
$b)-c)$ are replaced with, for example,
$\|f''_\S\|_{\rm HS} + |\nabla_\S f|\leq 1$ on the sphere.

\vskip5mm
{\bf Proof of Corollary 15.1.} It is rather similar to the proof of Corollary 11.2.
Replacing the spherical gradient with the operator $D$ defined in (13.4) and 
repeating the same arguments on the basis of the covariance identity (13.5), 
we obtain the analog of (11.4), namely
\be
\int_{\S^{n-1}} e^f\,d\sigma_{n-1} \leq \int_{\S^{n-1}}\int_{\S^{n-1}} 
\exp\{c_n \left<Df(x),D f(y)\right>\}\, d\mu_n(x,y)
\en
with constants $c_n$ satisfying (13.6). Moreover, using
\bee
\left<Df(x),D f(y)\right> 
 & \leq &
\|D f(x)\|_{\rm HS}\, \|D f(y)\|_{\rm HS} \\
 & \leq &
\frac{1}{2}\,\|D f(x)\|_{\rm HS}^2 + \frac{1}{2}\,\|D f(y)\|_{\rm HS}^2
\ene
and the property that the measure $\mu_n$ has $\sigma_{n-1}$ as marginals,
the double integral in (15.3) can be bounded by
\bee
\int_{\S^{n-1}} \int_{\S^{n-1}} \exp\Big\{\frac{c_n}{2}\,\|D f(x)\|_{\rm HS}^2 + 
\frac{c_n}{2}\,\|D f(y)\|_{\rm HS}^2\Big\}\, d\mu_n(x,y)
 & & \\
 & & \hskip-90mm \leq \
\Big(\int_{\S^{n-1}} \int_{\S^{n-1}} 
e^{c_n\,\|D f(x)\|_{\rm HS}^2}\, d\mu_n(x,y)\Big)^{1/2}
\Big(\int_{\S^{n-1}} \int_{\S^{n-1}} 
e^{c_n\,\|D f(y)\|_{\rm HS}^2}\, d\mu_n(x,y)\Big)^{1/2} \\
 & & \hskip-90mm = \
\int_{\S^{n-1}} e^{c_n\,\|D f\|_{\rm HS}^2}\, d\sigma_{n-1},
\ene
where we also applied Cauchy's inequality. It remains to note that, by (14.3),
$$
\|D f(x)\|_{\rm HS}^2 \, \leq \,
\Big(\|f_\S''(x)\|_{\rm HS} + 2\, |\nabla_\S f(x)|\Big)^2 \, \leq \,
2\,\|f_\S''(x)\|_{\rm HS} + 8\, |\nabla_\S f(x)|^2.
$$
\qed

\vskip5mm
\section{{\bf Covariance representations on the line}}
\setcounter{equation}{0}

\vskip2mm
\noindent
As we have already emphasized, covariance identities on the Euclidean
space $\R^n$ of dimension $n \geq 2$ with the usual gradient
exist for Gaussian measures, only. However, in dimension $n=1$ the situation is 
completely different. In fact, for any probability measure $\mu$ on the real line
a covariance identity such as
\be
{\rm cov}_\mu(u,v) = \int_{-\infty}^\infty \int_{-\infty}^\infty 
u'(x) v'(y)\,d\lambda(x,y)
\en
exists for a suitable measure $\lambda$ on the plane $\R \times \R$. In this and
next sections, we collect several results in this direction and refer 
an interested reader to \cite{B-D} for more details and historical references. 
Denote by $C^\infty_b$
the class of all functions $u:\R \rightarrow \R$ having $C^\infty$-smooth, 
compactly supported derivatives (in which case $u$ are bounded).

\vskip5mm
{\bf Proposition 16.1.} {\sl Given a probability measure $\mu$ 
on the real line, $(16.1)$ holds true for all $u,v \in C^\infty_b$ with a unique
locally finite measure $\lambda = \lambda_\mu$. This measure is non-negative and
absolutely continuous with respect to the Lebesgue measure on the plane.
}

\vskip5mm
Here ``locally finite" means that $\lambda$ is finite on compact sets
in the plane. In this case, both sides of (16.1) are well-defined and finite
for all functions in $C^\infty_b$. Moreover, (16.1) is extended to all locally 
absolutely continuous, complex-valued functions $u,v$ such that
$$
\int_{-\infty}^\infty \int_{-\infty}^\infty 
|u'(x)|\, |u'(y)|\,d\lambda(x,y) < \infty, \quad
\int_{-\infty}^\infty \int_{-\infty}^\infty 
|v'(x)|\, |v'(y)|\,d\lambda(x,y) < \infty,
$$
where the derivatives $u'$ and $v'$ are understood in the Radon-Nikodym sense. 
This condition guarantees that the double integral in (16.1) is finite and that
$u$ and $v$ belong to $L^2(\mu)$.

The equality (16.1) can be further generalized as a covariance identity
\be
\cov(u(X),v(Y)) = \int_{-\infty}^\infty \int_{-\infty}^\infty
u'(x) v'(y)\,H(x,y)\,dx\,dy
\en
for arbitrary random variables $X$ and $Y$, where
$$
H(x,y) = \P\{X \leq x, Y \leq y\} - \P\{X \leq x\}\, 
\P\{Y \leq y\}, \quad x,y \in \R.
$$
Here, the particular case of the identical functions $u(x)=x$ and $v(y)=y$ 
corresponds to the observation by H\"offding \cite{H} 
(provided that $X$ and $Y$ have finite second moments). 

Moreover, when $X=Y$, and the random variable $X$ is distributed according 
to $\mu$, (16.2) is reduced to (16.1). One may therefore conclude
that the mixing measure in (16.1) has density
\be
H_\mu(x,y) = \frac{d\lambda(x,y)}{dx\,dy} = 
F(x \wedge y)\,(1 - F(x \vee y)), \quad x,y \in \R,
\en
where
$$
F(x) = \P\{X \leq x\} = \mu((-\infty,x]), \quad x \in \R,
$$
is the distribution function associated to $\mu$,
with notations $x \wedge y = \min(x,y)$, $x \vee y = \max(x,y)$.

\vskip5mm
{\bf Definition 16.2.} Following \cite{B-D}, we call $\lambda = \lambda_\mu$ 
the H\"offding measure and its density $H = H_\mu$ the H\"offding kernel 
associated to $\mu$.

\vskip5mm
For example, if $\mu$ is the Bernoulli measure 
assigning the weights $p \in (0,1)$ and $q=1-p$ to the points $a < b$, then
$\lambda = pq U$ where $U$ is the uniform distribution on the square
$(a,b) \times (a,b)$.

Being applied with $u(x) = v(x) = x$, the identity (16.1) 
shows that the total mass of the H\"offding measure is the variance
$$
\lambda(\R \times \R) = \int_{-\infty}^\infty \int_{-\infty}^\infty 
H(x,y)\,dx\,dy = \Var(X).
$$
Thus, $\lambda$ is finite, if and only if $\mu$ has finite second moment.

Once the measure $\lambda$ is finite, it may also be described via its
Fourier-Stieltjes transform in terms of the characteristic function
of the random variable $X$,
$$
f(t) = \E\,e^{itX} = \int_{-\infty}^\infty e^{itx}\,d\mu(x),
\quad t \in \R.
$$
Namely, applying (16.1) to the exponential functions $u(x) = e^{itx}$ 
and $v(y) = e^{isy}$ ($t,s \in \R$), we obtain an explicit formula
\be
\hat \lambda(t,s) = \int_{-\infty}^\infty \int_{-\infty}^\infty 
e^{itx + isy}\,d\lambda(x,y) =
\frac{f(t) f(s) - f(t+s)}{ts} \quad (t,s \neq 0).
\en
In particular, this provides the uniqueness part in Theorem 1.1.

Thus, the expression on the above right-hand side represents a positive 
definite function in two real variables, as long as the characteristic 
function $f$ is twice differentiable. 

Moreover, there is a similar property of H\"offding's kernels: every such 
function $H = H_\mu$ defined in (16.3) is positive definite on the plane, that is,
$\sum_{i,j=1}^n a_i a_j H(x_i,x_j) \geq 0$ for any collection $a_i, x_i \in \R$. 
More generally,
$$
\int_{-\infty}^\infty \int_{-\infty}^\infty f(x) f(y) H(x,y)\,dx\,dy \geq 0,
$$
for any (bounded) measurable function $f$ on the real line.

Being positive definite, every H\"offding kernel satisfies
$H(x,y)^2 \leq H(x,x) H(y,y)$, which may be used to construct a pseudometric
$$
d(x,y) = \big(H(x,x) - 2H(x,y) + H(y,y)\big)^{1/2}.
$$
This property can be strengthened in terms of the H\"offding measure. 
In particular, 
$$
\lambda(A \times B)^2 \, \leq \, 
\lambda(A \times A) \, \lambda(B \times B)
$$
for all Borel sets $A,B \subset \R$.

Since the kernel $H = H_\mu$ is symmetric about the diagonal $x=y$, the
H\"offding measure $\lambda = \lambda_\mu$ has equal marginals 
$\Lambda = \Lambda_\mu$ defined by
\be
\Lambda(A) = \lambda(A \times \R) = \int_A^\infty \int_{-\infty}^\infty
H(x,y)\,dx\,dy, \quad A \subset \R \ ({\rm Borel}).
\en
Obviously, it is absolutely continuous with respect to the Lebesgue measure.
More precisely, we have the following density description. As before, assume
that the random variable $X$ is distributed according to $\mu$.

\vskip5mm
{\bf Proposition 16.3.} {\sl If $X$ has finite first absolute moment, then 
the marginal $\Lambda$ is finite and has density
\be
h(x) = \frac{d\Lambda(x)}{dx} = \int_x^\infty (y-a)\,dF(y), \quad
a = \E X.
\en
In particular, it is unimodal with mode at the point $a$, that is, $h(x)$ is 
non-decreasing on the half-axis $x < a$ and is non-increasing for $x > a$.
Moreover, it is continuous at $x=a$ with
\be
h(a-) = h(a+) = \frac{1}{2}\,\E\,|X-a|.
\en
}

\vskip2mm
In the case $\E\,|X| = \infty$, the density of $\Lambda$
is a.e. infinite on the support interval of $\mu$.

\vskip5mm
{\bf Proposition 16.4.} {\sl Assuming that $X$ has finite first absolute moment, 
the marginal $\Lambda$ is a multiple of $\mu$, if and only if $\mu$ is Gaussian.
}

\vskip5mm
Once $\Lambda = \sigma^2 \mu$ with some constant $\sigma^2$, the measures 
$\Lambda$ and $\lambda$ are necessarily finite, so that
$\mu$ must  have a finite second moment. If $\E X = 0$ (without loss of generality),
it follows from (16.4) that the Fourier-Stieltjes transform of 
$\Lambda$ is given by
$$
\hat \Lambda(t) = \hat \lambda(t,0) = 
\int_{-\infty}^\infty e^{itx}\,h(x)\,dx = -\frac{f'(t)}{t}, \quad 
t \in \R, \ t \neq 0,
$$
where $f$ is the characteristic function of $X$. Hence, $\Lambda = \sigma^2 \mu$ 
if and only if $f'(t) = -\sigma^2 t f(t)$, $t \in \R$. But this is only possible 
when $\mu$ is Gaussian with mean zero and variance $\sigma^2$.

Often, the marginals of H\"offding's measures appear in the particular case 
of the covariance representation (16.1) with the function $v(x) = x$. Then we arrive at
$$
\cov(X,u(X)) = \int_{-\infty}^\infty u'(x)\,h(x)\,dx,
$$
holding true as long as the integral is convergent.
If $\mu$ is supported on an interval $\Delta$ and has there an a.e.
positive density $p$, this formula may be rewritten as
\be
\cov(X,u(X)) = \E\,\tau(X) u'(X).
\en
Here the function
$$
\tau(x) = \frac{h(x)}{p(x)} =
\frac{1}{p(x)} \int_x^\infty (y-a)\,p(y)\,dy, \quad x \in \Delta,
$$
is called the Stein kernel. We have $\tau(x) = 1$ (a.e. on $\Delta$), if and 
only if $\mu$ is the standard Gaussian measure, in which case (16.8) becomes 
Stein's equation 
$$
\E\, X u(X) = \E\,u'(X). 
$$

After the pioneering work \cite{St}, 
the identity (16.8) served as a starting point for the extensive development 
of Stein's method as an approach to various forms of the central limit 
theorem and estimating the distances to the normal law.
For a comprehensive exposition of this theory, we refer an interested 
reader to the book \cite{C-G-S} and survey \cite{Sh}.

\vskip5mm
\section{{\bf Periodic covariance representations}}
\setcounter{equation}{0}

\vskip2mm
\noindent
The spherical covariance representation (1.4) of Theorem 1.3 in dimension
$n=2$, that is, on the circle $\S^1$, can be reduced to the covariance representation
\be
{\rm cov}_\mu(u,v) = \int_0^1 \int_0^1 u'(x) v'(y)\, d\lambda(x,y)
\en
for the uniform distribution $\mu = m$ on $(0,1)$ in the class of all
1-periodic smooth functions.

\vskip5mm
{\bf Definition 17.1.} We call a signed symmetric measure $\lambda$ 
on $[0,1) \times [0,1)$ a mixing measure for a given probability measure 
$\mu$ on $[0,1)$, if (17.1) holds true for all 1-periodic smooth functions 
$u$ and $v$ on the real line.

\vskip5mm
Following \cite{B-D}, let us now describe several results about this type
of covariance representations.
As we discussed before, the identity (17.1) always holds with the H\"offding measure
$\lambda = \lambda_\mu$. But, its marginals may be a multiple of $\mu$
in the Gaussian case only. This motivates the following question. Given $\mu$, 
is it possible to choose a mixing measure $\lambda$
whose marginals are multiples of $\mu$? If so, how to describe all of them
and choose a best one (in some sense)?

\vskip5mm
{\bf Proposition 17.2.} {\sl Let $\mu$ be a probability measure on $[0,1)$ with 
H\"offding's measure $\lambda_\mu$. Subject to the constraint that the 
marginal distribution of $\lambda$ in $(17.1)$ is equal to $c\mu$ for 
a prescribed value $c \in \R$, the mixing measure $\lambda$ exists, is unique, 
and is given by
\bee
\lambda 
 & = &
\lambda_\mu + (\sigma^2 - c)\,m \otimes m \nonumber \\
 & & + \ 
c\, (\mu \otimes m + m \otimes \mu) -
(\Lambda_\mu \otimes m + m \otimes \Lambda_\mu),
\ene
where $\Lambda_\mu$ is the marginal of $\lambda_\mu$ 
and $\sigma^2$ is the variance of $\mu$.
}

\vskip5mm
Here and in the sequel, $m$ denotes the uniform distribution on $(0,1)$.

One can specialize this characterization to the measue $m$
and consider identities
\be
{\rm cov}_m(u,v) = \int_0^1 \int_0^1 u'(x) v'(y)\, d\lambda(x,y).
\en
Using Proposition 16.1, Proposition 17.2 leads to the following assertion which
will be needed in the study of covariance identities on the circle.

\vskip5mm
{\bf Corollary 17.3.} {\sl Subject to the constraint that the marginal 
distribution of a mixing measure $\lambda$ in $(17.2)$ is equal 
to $cm$, $c \in \R$, the measure $\lambda$ is unique and has density
$$
\frac{\lambda(x,y)}{dx\,dy} = Q(|x-y|) + \Big(c - \frac{1}{24}\Big), 
\quad x,y \in (0,1),
$$
where
\be
Q(h) = \frac{1}{8}\,\big[\,1 - 4h(1 - h)\big], \quad 0 \leq h \leq 1.
\en
}

\vskip2mm
Note that $Q(h) \geq 0$ for all $h \in [0,1]$ and the inequality becomes
an equality for $h=\frac{1}{2}$.

Moreover, the mixing measure $\lambda$ is non-negative, if and only if 
$c \geq \frac{1}{24}$. Hence, the smallest positive measure (for the usual 
comparison) corresponds to the parameter $c = \frac{1}{24}$. In this sense, 
the optimal variant of (18.2) is given by the covariance representation
$$
{\rm cov}_{m}(u,v) = 
\frac{1}{24} \int_0^1 \int_0^1 u'(x) v'(y)\,d\nu(x,y)
$$
with a probability measure $d\nu(x,y) = 24\,Q(|x-y|)\,dx\,dy$
on $(0,1) \times (0,1)$.
It has the uniform distribution $m$ on $(0,1)$ as a marginal one.

If we want to write down a similar representation on 
the interval $(0,T)$, $T>0$, one may use a linear transform. 
Let $m_T$ denote the uniform distribution on $(0,T)$. Then we get that, 
for all smooth $T$-periodic functions $u$ and $v$,
and for all $c \geq 1/24$, 
\be
{\rm cov}_{m_T}(u,v) = \int_0^T \int_0^T u'(x) v'(y)\,d\lambda_T(x,y)
\en
with a positive measure having the density
\be
\frac{d\lambda_T(x,y)}{dx\,dy} = 
Q\Big(\Big|\frac{x}{T} - \frac{y}{T}\Big|\Big) + \Big(c - \frac{1}{24}\Big)
\en
on $(0,T) \times (0,T)$. It has the marginal $cT m_T(dx) = c\,dx$ on $(0,T)$.


\vskip10mm
\section{{\bf From the circle to the interval}}
\setcounter{equation}{0}

\vskip2mm
\noindent
Let us now return to the covariance representations on the circle
\be
{\rm cov}_{\sigma_1}(f,g) = \int_{\S^1}\int_{\S^1} 
\left<\nabla_\S f(x),\nabla_\S g(y)\right> d\nu(x,y),
\en
where we admit that $\nu$ may be a signed measure on the torus 
$\S^1 \times \S^1$. The question of whether such a measure $\nu$ exists
is settled in Proposition 6.6 with $\nu = \nu_2$.
Here we relate (18.1) to the similar covariance identity
\be
{\rm cov}_{m_{2\pi}}(u,v) = 
\int_0^{2\pi} \int_0^{2\pi} u'(t) v'(s)\,d\lambda(t,s)
\en
on the semi-open interval $[0,2\pi)$ with respect to the uniform
measure $m_{2\pi}$ in the class of all $2\pi$-periodic functions
$u$ and $v$ on the real line. 

Identifying $\R^2$ with the complex plane $\C$, to every smooth
function $f$ on $\S^1$ one associates
\be
u(t) = f(e^{it}) = f(\cos t,\sin t), \quad t \in \R,
\en
which is a smooth, $2\pi$-periodic function on the real line. 
Conversely, starting with such $u(t)$, (18.3) defines a smooth 
function $f(x)$ on the circle in a unique way. 
Given another smooth function $g$ on $\S^1$, define similarly
$v(t) = g(e^{it}) = g(\cos t,\sin t)$, $t \in \R$, so that
\be
{\rm cov}_{\sigma_1}(f,g) = {\rm cov}_{m_{2\pi}}(u,v).
\en

One can also rewrite the double integral in (18.1) explicitly in terms
of $u$ and $v$. Recall that the spherical gradient 
$w=\nabla_\S f(x)$ at the point $x \in \S^{n-1}$ for a smooth function 
$f$ represents the unique vector in $\R^n$ with the smallest 
Euclidean norm such that
$$
f(x') - f(x) = \left<w,x'-x\right> + o(|x'-x|), \quad x' \rightarrow x  \
(x' \in \S^{n-1}).
$$
In the case of the circle, one can put $x = e^{it}$, $x' = e^{i(t + \ep)}$
and rewrite the above as
\be
u(t+\ep) - u(t) = \left<w,e^{i(t+\ep)} - e^{it}\right> + o(\ep) =
\left<w,i e^{it}\right> \ep + o(\ep), \quad \ep \rightarrow 0.
\en
Since $i e^{it} = (-\sin t,\cos t)$, we have
$$
\left<w,i e^{it}\right> = -w_1 \sin t + w_2 \cos t, \quad w = (w_1,w_2).
$$ 
Thus, by (18.5),
\be
u'(t) = -w_1 \sin t + w_2 \cos t,
\en
implying that $|u'(t)| \leq |w|$. Moreover, the equality here is attained 
for $w_1 = -u'(t)\,\sin t$ and $w_2 = u'(t) \cos t$. These numbers are
therefore the coordinates of the shortest vector $w$ satisfying (18.5).
One may conclude that
$$
\nabla_\S f(x) = u'(t)\,(-\sin t,\cos t), \quad x = (\cos t, \sin t),
$$
and similarly
$$
\nabla_\S g(y) = v'(s)\,(-\sin s,\cos s), \quad y = (\cos s, \sin s),
$$
which gives
$$
\left<\nabla_\S f(x),\nabla_\S g(y)\right> = u'(t) v'(s)\,\cos(t-s).
$$

Denote by $\tilde \nu$ the measure on $[0,2\pi) \times [0,2\pi)$ 
such that $\nu$ is the image of $\tilde \nu$ under the map
$T(t,s) = (e^{it},e^{is})$. In view of (18.4), one may rewrite
(18.1) as a covariance identity
$$
{\rm cov}_{m_{2\pi}}(u,v) = 
\int_0^{2\pi} \int_0^{2\pi} u'(t) v'(s)\,\cos(t-s)\, d\tilde \nu(t,s),
$$
Thus, we obtain:

\vskip5mm
{\bf Lemma 18.1.} {\sl The representation $(18.1)$ on the circle with 
a signed measure $\nu$ on $\S^1 \times \S^1$ is equivalent to the covariance
representation $(18.2)$ in the class of all $2\pi$-periodic smooth functions 
$u$ and $v$ with a signed measure $\lambda$ on $[0,2\pi) \times [0,2\pi)$ defined by
\be
d\lambda(t,s) = \cos(t-s)\,d\tilde \nu(t,s), \quad t,s \in [0,2\pi).
\en
}

\vskip2mm
Note that $\cos(t-s) = 0$ if and only if $\left<x,y\right> = 0$ for
$x = e^{it}$, $y = e^{is}$. Hence, the map 
$\nu \sim \tilde \nu \rightarrow \lambda$ in (18.7)
transfers the part of the measure $\nu$ supported on the set 
\be
\Delta = \Big\{(x,y) \in \S^1 \times \S^1: \left<x,y\right> = 0\Big\}
\en
to zero. With this in mind, one may always assume 
that $\nu$ is supported outside $\Delta$.

To proceed, we will need to clarify the correspondence $\nu \sim \tilde \nu$
for an important class of measures on the torus.

\vskip5mm
{\bf Lemma 18.2.} {\sl A finite measure $\nu$ on the torus $\S^1 \times \S^1$ 
has density $\psi(x,y) = \psi(\left<x,y\right>)$ with respect to 
$\sigma_1 \otimes \sigma_1$ depending on the inner product $\left<x,y\right>$
if and only if $\tilde \nu$ has density 
\be
\frac{d\tilde \nu(t,s)}{dt\,ds} =
\frac{1}{(2\pi)^2}\,\psi(\cos(t-s)), \quad 0 < t,s < 2\pi,
\en
with respect to the Lebesgue measure on the square $[0,2\pi) \times [0,2\pi)$. 
}

\vskip5mm
{\bf Proof.} The assumption that $\nu$ appears as the image
of $\tilde \nu$ under the map $(t,s) \rightarrow (e^{it},e^{is})$
may be equivalently stated as an integral identity
\be
\int_{\S^1}\int_{\S^1} w(x,y)\,d\nu(x,y) =
\int_0^{2\pi} \int_0^{2\pi} w(e^{it},e^{is})\,d\tilde \nu(t,s),
\en
holding for all Borel measurable functions $w:\S^1 \times \S^1 \rightarrow \R$
(assuming that the integrals exist in the Lebesgue sense).
Since this map transfers the uniform distribution on $[0,2\pi) \times [0,2\pi)$
to the uniform distribution on the torus, we also have a similar general identity
\be
\int_{\S^1}\int_{\S^1} w(x,y)\,d\sigma_1(x)\,d\sigma_1(y) =
\frac{1}{(2\pi)^2} \int_0^{2\pi} \int_0^{2\pi} w(e^{it},e^{is})\,dt\,ds.
\en

If $\nu$ has density $\psi(x,y) = \psi(\left<x,y\right>)$ with respect 
to $\sigma_1 \otimes \sigma_1$, then (18.10) becomes
\be
\int_{\S^1}\int_{\S^1} w(x,y)\,\psi(\left<x,y\right>)\,
d\sigma_1(x)\,d\sigma_1(y) =
\int_0^{2\pi} \int_0^{2\pi} w(e^{it},e^{is})\,d\tilde \nu(t,s).
\en
On the other hand, applying (18.11) to functions of the form 
$w(x,y)\,\psi(\left<x,y\right>)$ and using $\left<x,y\right> = \cos(t-s)$ for 
$x = e^{it}$, $y = e^{is}$, we then obtain that
$$
\int_{\S^1}\int_{\S^1} w(x,y)\,\psi(\left<x,y\right>)\,
d\sigma_1(x)\,d\sigma_1(y) = \frac{1}{(2\pi)^2} \int_0^{2\pi} \int_0^{2\pi} 
w(e^{it},e^{is})\,\psi(\cos(t-s))\,dt\,ds.
$$
The right-hand side in this equality should be equalized with the 
right-hand side in (18.12). Since in these equalities $w$ may be 
an arbitrary bounded, Borel measurable function on the torus, the desired
relation (18.9) readily follows.

The argument in the opposite direction is similar.
\qed


\vskip5mm
\section{{\bf Mixing measures on the circle and the interval}}
\setcounter{equation}{0}

\vskip2mm
\noindent
How can one describe all measures $\nu$ in (18.1) under natural constraints? 
As we have already explained, $\nu$ is not unique: If a signed measure $\kappa$ 
is supported on the set $\Delta$ defined in (18.8), and (18.1) holds true with $\nu$, 
this relation will continue to hold for the measure $\nu + \kappa$. In addition, 
$\kappa$ may have $\sigma_1$ as marginals. However, being supported on the set 
of measure zero, any such measure $\kappa$ may not be absolutely continuous 
with respect to $\sigma_1 \otimes \sigma_1$.

Since (18.1) has been related to (18.2), the question of possible measures
$\nu$ may be reformulated in terms of the latter covariance representation 
in terms of $\lambda$. If we require that $\nu$ should have marginals being
multiples of $\sigma_1$, one may then apply Corollary 17.3.

Moreover, let us require that the measure $\nu$ in (18.1) has density 
$\psi(x,y) = \psi(\left<x,y\right>)$ with respect to $\sigma_1 \otimes \sigma_1$ 
depending on the inner product $\left<x,y\right>$ on the torus $\S^1 \times \S^1$. 
Recall that $\tilde \nu$ was defined 
as a signed measure on $[0,2\pi) \times [0,2\pi)$ such that $\nu$ appears as the 
distribution of $\tilde \nu$ under the map $(t,s) \rightarrow (e^{it},e^{is}$).
Then, by Lemma 18.2, 
$\tilde \nu$ has density $\frac{1}{(2\pi)^2}\,\psi(\cos(t-s))$ with respect 
to the Lebesgue measure on the square $[0,2\pi) \times [0,2\pi)$. Hence, 
by Lemma 18.1, the measure $\lambda$ has density 
$$
\frac{d\lambda(t,s)}{dt\,ds} = \frac{1}{(2\pi)^2}\,
\cos(t-s)\,\psi(\cos(t-s)), \quad t,s \in (0,2\pi).
$$
In particular, this measure is symmetric and has a marginal distribution 
proportional to $m_{2\pi}$. 

On the other hand, we are in position to apply Corollary 17.3 in the form
(17.5) with $T = 2\pi$ and conclude that, for some constant $c \in \R$,
$$
\frac{d\lambda(t,s)}{dt\,ds} = Q\Big(\frac{|t-s|}{2\pi}\Big) +c', 
$$
where 
$$
Q(h) = \frac{1}{8}\,\Big[\,1 - 4h(1 - h)\Big], \quad c' = c - \frac{1}{24}.
$$
In this case, $\lambda$ has as a marginal $c$ times the Lebesgue measure on 
$[0,2\pi)$. Equalizing the two formulas, we arrive at
\be
\frac{1}{(2\pi)^2}\,(\cos h)\,\psi(\cos h) = Q\Big(\frac{h}{2\pi}\Big) +c',
\en
which holds for almost all $h \in (0,2\pi)$, or equivalently a.e.
$$
\psi(\sin h) = \frac{(2\pi)^2}{\sin h}\,
\Big(Q\Big(\frac{h + \frac{\pi}{2}}{2\pi}\Big) +c'\Big), \quad
-\frac{\pi}{2} < h < \frac{3\pi}{2}.
$$
Here, $\sin h \sim h$ in the neighborhood of zero, over which
the left-hand side is integrable. However, as $h \rightarrow 0$,
$$
Q\Big(\frac{h + \frac{\pi}{2}}{2\pi}\Big) +c' \rightarrow
Q(1/4) + c' = \frac{1}{32} + c',
$$
which shows that, for the integrability of the right-hand side,
it is necessary that $c' = -\frac{1}{32}$.
Similarly, as $h \rightarrow \pi$,
$$
Q\Big(\frac{h + \frac{\pi}{2}}{2\pi}\Big) +c' \rightarrow
Q(3/4) + c' = \frac{1}{32} + c',
$$
so that we do not obtain a new restriction. Note that 
\bee
K(h) \, \equiv \, Q(h) - \frac{1}{32} 
 & = &
\frac{1}{8}\,\Big[\,1 - 4h(1 - h)\Big] - \frac{1}{32} \\
 & = &
\frac{3}{32} - \frac{1}{2}\,h(1 - h) \ = \
\frac{1}{2}\,\Big(h - \frac{1}{4}\Big) \,\Big(h - \frac{3}{4}\Big).
\ene
This new kernel is symmetric about the point $1/2$, with
$$
K(0) = K(1) = \frac{3}{32}, \quad K(1/2) = - \frac{1}{32}, \quad
K(1/4) = K(3/4) = 0,
$$ 
so that it is positive in 
$0 < h < 1/4$ and $3/4 < h < 1$, and is negative in the interval 
$1/4 < h < 3/4$. This behavior is rather similar
to the one of the function $\cos(\frac{h}{2\pi})$.

Thus, with the necessary value $c' = -\frac{1}{32}$ from (19.1)
we obtain that
\be
\psi(\cos h) = (2\pi)^2\, \frac{K\big(\frac{h}{2\pi}\big)}{\cos h},
\quad 0 \leq h \leq 2\pi,
\en
In particular, applying this equality with $h=\pi$, $h = \frac{\pi}{2}$, 
and $h = 0$, we have
$$
\psi(-1) = \frac{\pi^2}{8}, \quad 
\psi(0) = \frac{\pi}{2}, \quad \psi(1) = \frac{3\pi^2}{8}.
$$
In this case, 
$c = c' + \frac{1}{24} = -\frac{1}{32} + \frac{1}{24} = \frac{1}{96}$
which is the total mass of the measure $\lambda$.

Let us now describe the total mass of the measure $\nu$, that is,
$\tilde \nu$. Applying (18.12) with $w=1$, we have
\begin{eqnarray}
\nu(\S^1 \times \S^1)
 & = &
\int_{\S^1}\int_{\S^1} \psi(\left<x,y\right>)\,d\sigma_1(x)\,d\sigma_1(y) \nonumber \\
 & = &
\frac{1}{(2\pi)^2} \int_0^{2\pi}\! \int_0^{2\pi} \psi(\cos(t-s))\,dt\,ds \, = \, 
\int_0^1\! \int_0^1 \psi(\cos(2\pi(t-s)))\,dt\,ds \nonumber \\
 & = &
\int_{-1}^1 \psi(\cos(2\pi h))\,(1-|h|)\,dh \, = \,
2 \int_0^1 \psi(\cos(2\pi h))\,(1-h)\,dh \nonumber \\
 & = &
2\,(2\pi)^2 \int_0^1 \frac{K(h)}{\cos(2\pi h)}\,(1-h)\,dh.
\end{eqnarray}
One can now refine Proposition 6.1, in which necessarily $\psi_2 = \psi$.

\vskip5mm
{\bf Proposition 19.1.} {\sl On the torus $\S^1 \times \S^1$ there exists 
a unique measure $\nu$ with density of the form $\psi(\left<x,y\right>)$ with 
respect to $\sigma_1 \otimes \sigma_1$ such that, for all smooth functions 
$f,g$ on $\S^1$,
\be
{\rm cov}_{\sigma_1}(f,g) = \int_{\S^1}\int_{\S^1} 
\left<\nabla_\S f(x),\nabla_\S g(y)\right> d\nu(x,y).
\en
The function $\psi(\alpha)$ is positive, increasing,
continuous for $|\alpha| \leq 1$, and is given by $(19.2)$ with
$$
K(h) = \frac{1}{32}\,(4h-1)(4h-3), \quad 0 \leq h \leq 1.
$$
In addition $\frac{\pi^2}{8} \leq \psi(\alpha) \leq \frac{3\pi^2}{8}$, where
equalities are attained for $\alpha=-1$ and $\alpha=1$.
}

\vskip5mm
One may rewrite (19.4) equivalently as
$$
{\rm cov}_{\sigma_1}(f,g) = c \int_{\S^1}\int_{\S^1} 
\left<\nabla_S f(x),\nabla_S g(y)\right> d\mu(x,y)
$$
with the assumption that $\mu$ is a probability measure on $\S^1 \times \S^1$
with marginals equal to $\sigma_1$. Here the constant $c$ is described in
(19.3).

As an alternative variant for (19.5), we have the following modified
representation.

\vskip5mm
{\bf Proposition 19.2.} {\sl On the torus $\S^1 \times \S^1$ there exists a finite 
positive measure $\nu$ such that, for all smooth functions  $f,g$ on $S^1$,
\be
{\rm cov}_{\sigma_1}(f,g) = \int_{\S^1}\int_{\S^1} 
\frac{\left<\nabla_\S f(x),\nabla_\S g(y)\right>}{\left<x,y\right>}\,
d\nu(x,y).
\en
}

\vskip5mm
{\bf Proof.} Recall that in the case $\left<x,y\right> = 0$, necessarily
$\left<\nabla_\S f(x),\nabla_\S g(y)\right> = 0$. Hence, there is 
no uncertainty in the integrand in (19.5).
Moreover, for $x = e^{it}$, $y = e^{is}$, we have
$\left<x,y\right> = \cos(t-s)$. Transferring the circle to $[0,2\pi)$ 
and the torus to the square $[0,2\pi) \times [0,2\pi)$ via 
the inverse of the map $T$, the representation (19.5) will take the form
$$
{\rm cov}_{m_{2\pi}}(u,v) = \int_0^{2\pi} \int_0^{2\pi} u'(t) v'(s)\,
d\tilde \nu(t,s)
$$
for the functions $u(t) = f(e^{it})$ and $v(s) = g(e^{is})$.

Here for $\tilde \nu$ one may take the H\"offding measure 
on the square $(0,2\pi) \times (0,2\pi)$ in the covariance representation
for the uniform measure $m_{2\pi}$, that is, with positive density
$$
\frac{d\tilde \nu(t,s)}{dt\,ds} = \psi(t,s) = F(t \wedge s)\,(1 - F(t \vee s)), \quad
t,s \in (0,2\pi),
$$
where $F(t) = t/(2\pi)$ is the distribution function for $m_{2\pi}$.
The measure $\nu$ in (19.5) will then appear as the image of $\tilde \nu$ 
under the map $T$.

\vskip2mm
{\bf Remark 19.3.}
As was mentioned in section 16, the measure $\nu$ has a total mass
$$
\nu(\S^1 \times \S^1) = \tilde \nu\big((0,2\pi) \times (0,2\pi)\big) =
\int_0^{2\pi} \int_0^{2\pi} \psi(t,s)\,dt\,ds = \Var(U) = \frac{\pi^2}{3},
$$
where $U$ is a random variable distributed according to $m_{2\pi}$. 
The marginal distributions of $\tilde \nu$
coincide and have a positive density on $(0,2\pi)$ given by
$$
\int_0^{2\pi} \psi(t,s)\,ds = \frac{1}{4\pi}\,t(2\pi - t), \quad
0 < t < 2\pi.
$$
After transferring this marginal distribution to the circle, we will
obtain the marginal distribution of $\nu$ which is however not
the uniform measure $\sigma_1$.

Note also that $\psi$ is vanishing on the boundary of the square
$[0,2\pi] \times [0,2\pi]$. Hence $\nu$ is absolutely continuous and
has a continuous density on the torus with respect to the product
measure $\sigma_1 \otimes \sigma_1$. However, in contrast with higher
dimensions, this density does not represent a function of the inner product
$\left<x,y\right>$ (since $\psi(t,s)$ is not a function of $\cos(t-s)$).

\vskip5mm
{\bf Acknowledgement.} Research of the first author was partially supported
by NSF grant DMS-2154001.

\end{document}